\documentclass[11pt]{article}

\usepackage{article_style}
\usepackage{writing}
\usepackage{adjustbox}
\usepackage{pdflscape}

\nocite{MAGMA}

    \makeatletter
    \@ifundefined{origwidetilde}{%
      \let\origwidetilde\widetilde
      \newcommand{\lw@drop}{0.2ex}%
      \newcommand{\lowtilde@}[2]{%
        \setbox\z@\hbox{$\m@th#1#2$}
        \setbox\tw@\hbox{$\m@th#1\origwidetilde{\vphantom{#2}\kern\wd\z@}$}
        \rlap{\raise-\lw@drop\box\tw@}
        \box\z@
      }%
      \renewcommand{\widetilde}[1]{\mathord{\mathpalette\lowtilde@{#1}}}%
    }{}%
    \makeatother

\title{Random homomorphisms from affine to spherical reflection groups} 

\author{Isabelle Steinmann}

\date{}
\begin{document}

\maketitle

\begin{abstract}
    The goal of this paper is to understand the set $\Hom(\widetilde{W},W)$ of homomorphisms from an affine reflection group $\widetilde{W}$ to its associated spherical reflection group $W$. We do this in three ways: numerically, by deriving an explicit formula for $|\Hom(\widetilde{W},W)|$; probabilistically, by exploring the question \textit{what does a random homomorphism $\widetilde{W}\to W$ look like?}; and algebraically, by explicitly writing down all possible homomorphisms on the standard generators.
    For example, we show that as $n\to\infty$ the probability that a random homomorphism $\widetilde{W}_n\to W_n$ is surjective tends to $\frac{3}{4}$ if $W_n$ is of type $C_{2n}$ or $D_n$, to $\frac{9}{16}$ if $W_n$ is of type $C_{2n+1}$, and to $1$ if $W_n$ is of type $A_n.$
\end{abstract}

\section{Introduction}
\label{sec: intro}

    The affine reflection groups $\widetilde{W}_n$ were classified by Coxeter {\rm\cite{Coxeter_affine}}. They consist of four infinite families of groups, $\widetilde{A}_n$, $\widetilde{B}_n$, $\widetilde{C}_n$, and $\widetilde{D}_n$, along with five exceptional groups, $\widetilde{E}_6$, $\widetilde{E}_7$, $\widetilde{E}_8$, $\widetilde{F}_4$ and $\widetilde{G}_2$. 
    The coroot lattice $Q(R^{\vee})\leqslant \widetilde{W}_n$ of translations in $\widetilde{W}_n$ gives rise to the decomposition
    $ \widetilde{W}_n=Q(R^{\vee})\rtimes W_n$, and hence to the natural projection 
    \begin{align*}
    \label{equ: standard proj}
        \pi_{\mathrm{standard}}:\widetilde{W}_n\cong Q(R^{\vee})\rtimes W_n&\to W_n \\
        (v,w)&\mapsto w.
    \end{align*}
    In this paper, we describe the homomorphisms $\widetilde{W}_n\to W_n$ by computing $|\Hom(\widetilde{W}_n,W_n)|$, and by examining a random homomorphism $\varphi\in \Hom(\widetilde{W}_n,W_n).$

\noindent{\bf Counting homomorphisms.} For any affine reflection group $\widetilde{W}_n$, there are finitely many homomorphisms $\widetilde{W}_n\to W_n$. Depending on the type of $\widetilde{W}_n$, there are many surjective homomorphisms $\widetilde{W}_n\to W_n$ other than $\pi_{\mathrm{standard}}$, even up to postcomposition with an automorphism of $W_n$.

    \begin{theorem} 
    \label{thm: main theorem}
    The number of homomorphisms $|\Hom(\widetilde{W}_n,W_n)|$ is:
    \begin{enumerate}[label=(\alph*)]
           
        \item For $A_n\cong \FS_{n+1}$, the order-$(n+1)!$ symmetry group of the $n$-simplex, when $n\geq4$ and $n\neq 5$:
        $$|\Hom(\widetilde{A}_n,A_n)|=(n+1)!+\sum_{k=0}^{\lfloor\frac{n+1}{2}\rfloor} 
        \frac{(n+1)!}{2^k\ k!\ (n+1-2k)!}.$$
        
        \item  For $B_n\cong  C_n\cong \FC_2^n\rtimes \FS_n$, the order-$2^nn!$ symmetry group of the $n$-cube, when $n\geq3$ and $n\neq 4,6$:
            \begin{align*}|\Hom(\widetilde{B}_n,B_n)|&= 2^{n+3}n!+        
            \sum_{k=0}^{\lfloor \frac{n}{4}\rfloor} \;
            \sum_{l=0}^{\lfloor \frac{n}{2}\rfloor-2k} 
            \frac{(\frac{3}{8})^l2^{2n}n!}{2^{7k} \ k!\ l!\ (n-4k-2l)!}  \\ 
            |\Hom(\widetilde{C}_n,C_n)|&= 2^{n+4}n!+
            \sum_{k=0}^{\lfloor \frac{n}{8}\rfloor} \;
            \sum_{l=0}^{\lfloor \frac{n}{4}\rfloor-2k}\; 
            \sum_{m=0}^{\lfloor \frac{n}{2}\rfloor-4k-2l} 
            \frac{2^{3n}7^{l+m}n!}{2^{20k+10l+4m} \ k!\ l!\ m!\ (n-8k-4l-2m)!}     . \end{align*}
        \item For $D_n\cong \FC_2^{n-1}\rtimes \FS_n$, the order-$2^{n-1}n!$ index $2$ subgroup of $C_n$ consisting of elements which contain an even number of coordinate inversions, when $n\geq3$ and $n\neq 4,6:$
        $$|\Hom(\widetilde{D}_n,D_n)|=2^{n+1+\frac{1+(-1)^n}{2}}n!+
        \frac{n!}{\lfloor\frac{n}{2}\rfloor!}+\sum_{k=0}^{\lfloor \frac{n}{2}\rfloor-1} \frac{2^{n-1}n!}{2^{2k}\ k!\ (n-2k)!}.$$
        \item For the exceptional affine groups, and the members of the infinite families omitted above:
        \begin{table}[h]
            \centering
            \begin{tabular}{l l l}
                $|\Hom(\widetilde{A}_2,A_2)|=28$ & $|\Hom(\widetilde{A}_3,A_3)|=106$ & $|\Hom(\widetilde{A}_5,A_5)|=1516$ \\ $|\Hom(\widetilde{B}_2,B_2)|=168$ & $|\Hom(\widetilde{B}_4,B_4)|=14176$ & $|\Hom(\widetilde{B}_6,B_6)|=660736$  \\ 
                $|\Hom(\widetilde{C}_4,C_4)|=54112$ & $|\Hom(\widetilde{C}_6,C_6)|=17499904$&  $|\Hom(\widetilde{D}_4,D_4)|=12332$ \\ 
                $|\Hom(\widetilde{D}_6,D_6)|=231152$&  $|\Hom(\widetilde{E}_6,E_6)|=52732$ &  $|\Hom(\widetilde{E}_7,E_7)|=2913248$\\
                 $|\Hom(\widetilde{E}_8,E_8)|=696929552$ & $|\Hom(\widetilde{F}_4,F_4)|=40096$ & $|\Hom(\widetilde{G}_2,G_2)|=88$
            \end{tabular}
        \end{table}
    \end{enumerate}
    \end{theorem}

\noindent{\bf Random homomorphisms.} The proof of Theorem \ref{thm: main theorem}  does more than compute $|\Hom(\widetilde{W}_n,W_n)|$: it also describes the homomorphisms $\widetilde{W}_n\to W_n$ in terms of their kernels and images. Hence, we can study properties of these homomorphisms to answer the question:
    
    \leftskip1.5cm\relax
    \rightskip1.5cm\relax
    \begin{center}
    \textit{What does a random homomorphism $\widetilde{W}_n\to W_n$ look like?} \\ 
    \end{center}
    \leftskip0cm\relax
    \rightskip0cm\relax
    
    \begin{theorem}
        \label{thm: expected image order A,B,C,D}
        Suppose $W_n$ is $A_n$, $B_n$, $C_n$, or $D_n$. Most of the homomorphisms $\widetilde{W}_n\to W_n$ have large images: asymptotically, the expected order of the image of a random $\varphi\in \Hom(\widetilde{W}_n,W_n)$ is:
        $$  \EE_{\varphi\in \Hom(\widetilde{W}_n,W_n)} \big[\ |\im(\varphi)|\  \big]\sim \begin{cases}
            |W_n| & W_n=A_n,\\
            \frac{3}{4}|W_n| & W_n=C_{2n},D_n,\\
            \frac{21}{32}|W_n| & W_n=C_{2n+1},\\
            \frac{5}{8}|W_n| & W_n=B_{2n},\\
            \frac{9}{16}|W_n| & W_n=B_{2n+1}.\\
        \end{cases} $$
        Moreover, the standard deviation of the order of the image is asymptotically
        $$  \sigma_{\varphi\in \Hom(\widetilde{W}_n,W_n)} \big[\  |\im(\varphi)|\  \big]\sim \begin{cases}
            0 & W_n=A_n,\\
            \frac{\sqrt{3}}{4}|W_n| & W_n=C_{2n},D_n,\\
            \frac{\sqrt{183}}{32}|W_n| & W_n=C_{2n+1},\\
            \frac{\sqrt{11}}{8}|W_n| & W_n=B_{2n},\\
            \frac{\sqrt{39}}{16}|W_n| & W_n=B_{2n+1}.\\
        \end{cases} $$
    \end{theorem}

    In addition to the expected image-order of a random homomorphism being large, it is also likely that a random homomorphism is surjective and hence has maximum image-order.
    
    \begin{theorem}
        \label{thm: prob surjective hom ABCD}
        Suppose $W_n$ is $A_n$, $B_n$, $C_n$, or $D_{n}$. Then 
        \begin{align*}
        \lim_{n\to\infty} \PP_{\varphi\in \Hom(\widetilde{W}_n,W_n)} \big[ \text{$\varphi$ is surjective} \big]=\begin{cases}
            1 & W_n=A_n,\\
            \frac{3}{4} & W_n=C_{2n},D_n,\\
            \frac{9}{16} & W_n=C_{2n+1},\\
            \frac{1}{2} & W_n=B_{2n}, \\ 
            \frac{3}{8} & W_n=B_{2n+1} .\\
        \end{cases}
        \end{align*}
    \end{theorem}

    Figure \ref{fig: proportional Hom(Cn~,Cn)} exhibits trends in the order of the image of  random homomorphism $\widetilde{C}_n\to C_n$.     
    \begin{figure}[h]
        \centering
        \includegraphics[width=\linewidth]{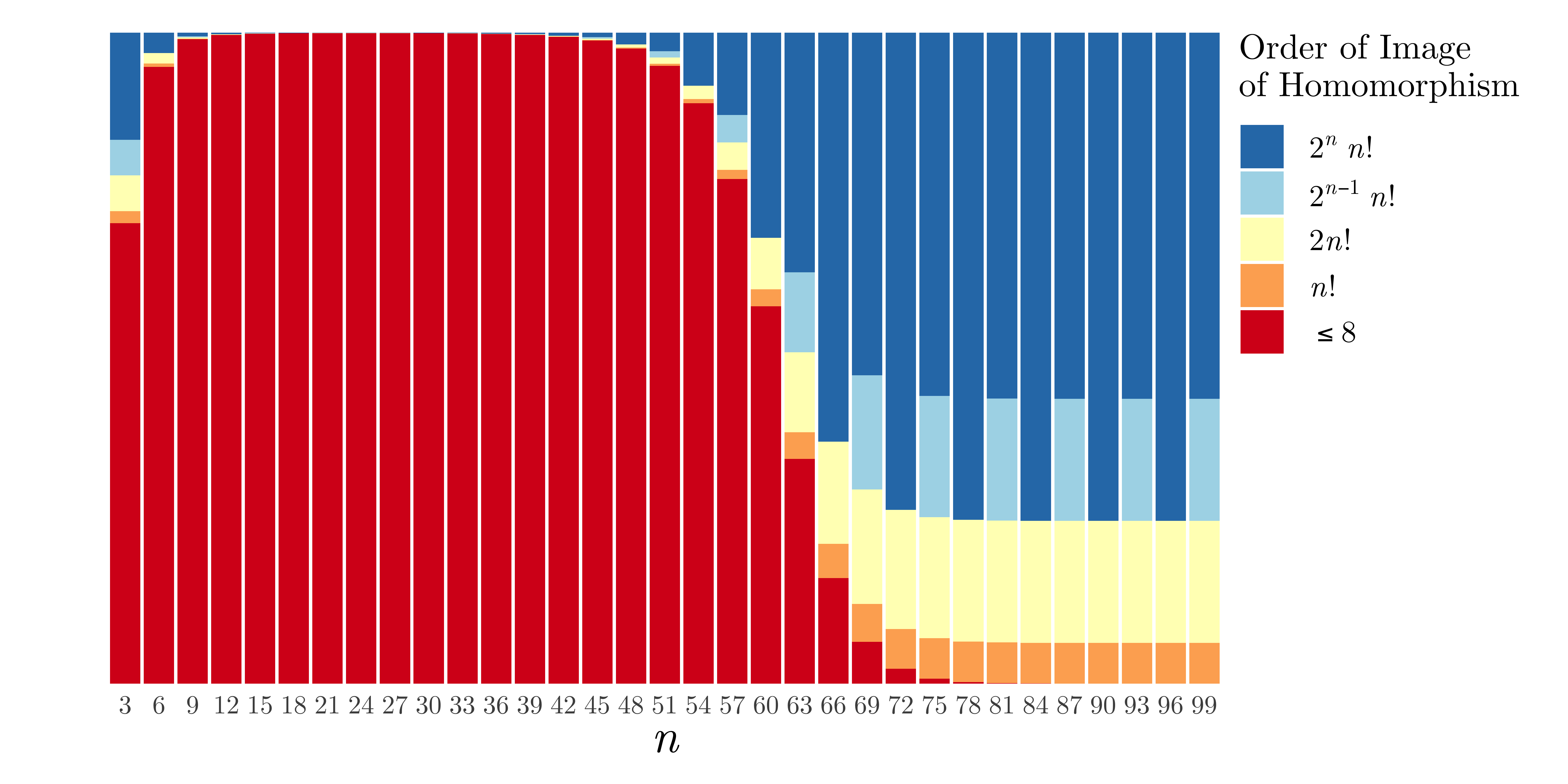}
        \caption{\small{A stacked bar graph which shows the proportional orders of the images of homomorphisms $\widetilde{C}_{n}\to C_{n}$ for $n=3m$ where $1\leq m \leq 33$.
        For example, when $n=3$ the graph shows that roughly two-thirds of the homomorphisms $\widetilde{C}_3\to C_3$ have image with order at most $8$. 
        We see that for a random $\varphi\in \Hom(\widetilde{C}_n,C_n)$ the probability that $|\varphi(C_n)|\leq8$ tends to $0$ as $n\to \infty$, but that this behaviour isn't evident until $n\approx 70$.
        Note that for $\varphi\in \Hom(\widetilde{C}_n,C_n)$, 
        $\varphi$ is surjective if and only if $|\im(\varphi)|=2^nn!=|C_n|$. Thus eventually $\frac{3}{4}$ of $\varphi\in \Hom(\widetilde{C}_{2n},C_{2n})$ are surjective, 
        and $\frac{9}{16}$ of  
        $\varphi\in \Hom(\widetilde{C}_{2n+1},C_{2n+1})$ are surjective.}}
        \label{fig: proportional Hom(Cn~,Cn)}
    \end{figure}

    A natural question is how generic $\pi_{\mathrm{standard}}$ is in the set of all surjective homomorphisms. Aside from a few exceptional cases when $n=2,4,6$, the proportion of surjective homomorphisms which factor through $\pi_{\mathrm{standard}}$ does not depend on $n$ and hence this probability can be stated directly, not merely asymptotically.

    \begin{theorem}
    \label{thm: prob surjective factors thru standard ABCD}
         Suppose $W_n$ is $A_n$, $B_n$, $C_n$, or $D_{n}$ and $n\neq 2,4,6$. Then 
        \begin{align*}
        \PP_{\varphi\in \Hom(\widetilde{W}_n,W_n)} \big[ \text{$\varphi$ factors through $\pi_{\mathrm{standard}}$ $\mid$ $\varphi$ is an surjective} \big]=\begin{cases}
            1 & W_n=A_n,\\
            \frac{1}{2} & W_n=B_{2n},\\
            \frac{1}{3} & W_n=B_{2n+1},D_n,\\
            \frac{1}{6} & W_n=C_{2n}, \\ 
            \frac{1}{9} & W_n=C_{2n+1} .\\
        \end{cases}
        \end{align*}
    \end{theorem}

     More generally, for any automorphism $f\in \Aut(W_n)$, the composition $f\circ \pi_{\mathrm{standard}}$ lies in $\Hom(\widetilde{W}_n,W_n)$. Hence we can also ask how generic these maps which factor through $\pi_{\mathrm{standard}}$ are in $\Hom(\widetilde{W}_n,W_n)$. 
     
    \begin{theorem}
    \label{thm: prob factors thru standard ABCD}
         Suppose $W_n$ is $A_n$, $B_n$, $C_n$, or $D_{n}$. Then 
        \begin{align*}
        \lim_{n\to \infty}\PP_{\varphi\in \Hom(\widetilde{W}_n,W_n)} \big[ \text{$\varphi$ factors through $\pi_{\mathrm{standard}}$} \big]=\begin{cases}
            1 & W_n=A_n, \\
            \frac{1}{2} & W_n=D_n,B_n,\\
            \frac{1}{4} & W_n=C_n.\\
        \end{cases}
        \end{align*}
    \end{theorem}

    \begin{remark}[Explicit description of homomorphisms]
        As well as these probabilistic descriptions, in Section \ref{sec: explicit homomorphisms} we also explicitly state homomorphisms $\widetilde{W}_n\to W_n$ on the standard generating set of $\widetilde{W}_n$ and, for types $B_n$, $C_n$, and $D_n$, on an arbitrary element of $\widetilde{W}_n$. For types $A_n$, $E_6$, $E_7$, and $E_8$, every homomorphism factors through $\pi_{\mathrm{standard}}$ and these are described in Theorem \ref{thm: explicit homs types A and E}. For types $F_4$ and $G_2$, there exist homomorphisms which do not factor through $\pi_{\mathrm{standard}}$ and we write down all homomorphisms on standard generators. For types $B_n$, $C_n$ and $D_n$, $\Hom(\widetilde{W}_n,W_n)$ is complex and there are many homomorphisms which do not factor through $\pi_{\mathrm{standard}}$. In Theorem \ref{thm: explicity homs types BCD}, we construct simple homomorphisms on subgroups of $\widetilde{W}_n$ which we combine to obtain homomorphisms on $\widetilde{W}_n$ itself. 
    \end{remark}

    \noindent
    {\bf Acknowledgements.}
    I would like to thank Benson Farb for his suggestion of this problem and for his constant guidance, encouragement, and comments. I would also like to thank Tai Lohrer for helping to create Figure \ref{fig: proportional Hom(Cn~,Cn)}.

    \noindent
    {\bf Computational resource disclosure.} Claude Opus 4.8 was used for proofreading, understanding the technical details of {\rm\cite{maxwell_normal_subs}}, checking the proofs in Section \ref{sec: explicit homomorphisms}, and in developing the proof of Lemma \ref{lem: exceptional quotients of Cn~}, which led to some of the author's ideas for Section \ref{sec: explicit homomorphisms}. AI was not used for any writing, for any proofs in Sections \ref{sec: hom tables ABCD}-\ref{sec: random homs}, or for the creation of any of the tables. 

\section{Preliminary material}

    In this section we introduce notation and background material, and discuss the structure of the proof of Theorem \ref{thm: main theorem}.

    Suppose $\widetilde{W}_n$ is an affine reflection group and $W_n$ is its associated spherical reflection group. We follow the standard notation of Bourbaki, as in Section 4 of {\rm\cite{bourbaki_groupes_1981}}, for spherical and affine reflection groups. This background information can also be found in Chapters 1,2 and 4 of {\rm\cite{humphreys_reflection_grps}}, although Humphreys uses slightly different notation.
    We denote by $Q(R^{\vee})$ the coroot lattice of $W_n$ and, by abuse of notation, we will also use $Q(R^{\vee})$ to signify the rank $n$ abelian subgroup of $\widetilde{W}_n$ consisting of translations by $v\in Q(R^{\vee})\cong \ZZ^n$. Note that $W_n$ acts on $Q(R^{\vee})$, and moreover
    $$ \widetilde{W}_n=Q(R^{\vee})\rtimes W_n. $$
    Let $\pi_{\text{standard}}:\widetilde{W}_n\to W_n$ be the projection map with $\pi_{\mathrm{standard}}(v,w)=w$.
    A simple root system of $W_n$ contains a unique longest root which we denote $\widetilde{\alpha}$. Then, if $s_1,\dots,s_n$ is the standard generating set of $W_n$, $\widetilde{W}_n$ is generated by $s_0,s_1,\dots,s_n$ where $s_0$ is the product of a translation by $\frac{2}{\lb \widetilde{\al},\widetilde{\al}\rb}\widetilde{\al}$ and $r_{\widetilde{\alpha}}$. On the generating set $s_0,\dots,s_n$ of $\widetilde{W}_n$, $\pi_{\mathrm{standard}}$ is: 
    $$ \pi_{\mathrm{standard}}(s_0)=r_{\widetilde{\alpha}}\qquad \pi_{\mathrm{standard}}(s_i)=s_i \quad \text{for }i\in\{1,\dots,n\} . $$

    Every reflection group $R$ has an index-2 \emph{rotation subgroup} (or \emph{even subgroup}) $R^+\unlhd R$ consisting of elements in $R$ that can be written in an even number of generating reflections. Hence for any group $G$ and $g\in G$ with $g^2=1$, we can define 
    \begin{equation}
    \label{equ: sgn:R-><g>}
        \mathrm{sgn}_R^g:R\to \lb g\rb\leqslant G
    \end{equation}
    to be the unique homomorphism onto $\lb g\rb$ with kernel contains $R^+$ (and is equal to $R^+$ if $g\neq 1$). The centre of any spherical reflection group $W_n$ has order at most $2$, so let $z_0$ be the (possibly trivial) generator of $Z(W_n)$. When $z_0\neq1$, we also write $z_0=-1$. Note that $\mathrm{sgn}_{W_n}^{z_0}$ is an endomorphism of $W_n$. 
    
\subsection{Method}

    There is a bijection between $\Hom(G,H)$ and the set of tuples $(K,M,\varphi)$, where $K\unlhd G$ is normal, $M\leqslant H$ is isomorphic to $G/K$, and $\varphi$ is an automorphism of $G/K$. 
    Note that normal subgroups $K\unlhd G$ are merely \emph{potential} kernels because, if there is no subgroup of $H$ isomorphic to $G/K$, then there is no homomorphism from $G\to H$ with kernel $K$. The normal subgroups of spherical and affine reflection groups and the quotient groups they induce are described by Maxwell in {\rm\cite{maxwell_normal_subs}}, and the automorphisms of spherical and affine reflection groups are well-known and are described in {\rm\cite[Ch. 2]{AutW_Franzsen}}.

    The number of subgroups of $H$ isomorphic to $\FC_2^k$ can be computed using $\V_1(H),\dots,\V_k(H)$, where 
    $$\V_k(H)=\big|\{(x_1,\dots,x_k)\in H^k: x_i^2=[x_i,x_j]=1\; \text{for }i,j\in \{1,\dots,k\}\} \big|$$ 
    For example, suppose $K\unlhd G$ is an index 2 subgroup (e.g. the rotation subgroup of a reflection group). The set of homomorphisms $G\to H$ with kernel $K$ is in bijection with the set of involutions of $H$, and the number of involutions in $H$ is equal to $V_1(H)-1$. 
    More generally, if $H$ is finite then we can compute the order of $\Hom(G,H)$ by:
    \begin{equation}
        \big|\mathrm{Hom}(G,H)\big|=\sum_{\substack{ K\unlhd G\\ |G:K|\leq |H|}} |\{X\leqslant H: X\cong G/K\}|\cdot |\Aut(G/K)|.
    \end{equation}

    \begin{definition}
    \label{def: homomorphism table}
    For any pair of groups $(G,H)$, one can construct a \emph{homomorphism table} which collates data on the homomorphisms $G\to H$. The rows of this table correspond to normal subgroups $K\unlhd G$ with $|G:K|\leq |H|$. The first two columns state $|G:K|$ and give a description of $K$. The third column describes the isomorphism type of the quotient group $G/K$, and hence the isomorphism type of the image of a homomorphism from $G\to H$ with kernel $K$. 
    The fourth gives the number of subgroups of $H$ isomorphic to $G/K$. These are (setwise) images of homomorphisms from $G$ to $H$ with kernel $K$. The fifth lists the order of the automorphism group of $G/K$. This counts the number of homomorphisms from $G$ to $H$ with kernel $K$ and (setwise) fixed image. The final column counts the total number of homomorphisms from $G$ to $H$ with kernel $K$, and it is obtained by multiplying the entries of the previous two columns. 
    \end{definition}
    
    The homomorphism table of $(G,H)$ not only allows us to compute $|\Hom(G,H)|$ by summing the entries in the final column, but it also gives a classification of homomorphisms from $G$ to $H$ in terms of their kernels. 
    
\subsection{Normal subgroups of affine reflection groups}
\label{subsec: normal subs of W~}

    We aim to classify homomorphisms $\varphi:\widetilde{W}\to W$ in terms of their kernels, images, and the automorphisms of their images. Hence we are interested in normal subgroups $N\unlhd \widetilde{W}$ with $|\widetilde{W}:N|\leq |W|$, and the isomorphism types of $\widetilde{W}/N$, both of which are provided by Maxwell  in Tables 5, 6, and 7 of {\rm\cite{maxwell_normal_subs}} (except for the isomorphism type of $\widetilde{W}/N$ in a few cases, which are discussed in Lemma \ref{lem: exceptional quotients of Cn~}). Theorem \ref{thm: normal subs of W~} describes the sorts of normal subgroups that affine reflection group can have. Suppose $(W,S)$ and $(W',S')$ are Weyl groups. A homomorphism $\psi:W\to W'$ is a \emph{graph homomorphism} if, for every $s\in S$, $\psi(s)=s'$ or $\psi(s)=1$; and for every $s'\in S'$, there is some $s\in S$ such that $\psi(s)=s'$. Graph homomorphisms are surjective, and are evident from the Coxeter graphs of $W$ and $W'$.

    \begin{theorem}{\rm\cite[p. 360]{maxwell_normal_subs}}
        Suppose $\widetilde{W}$ is an affine reflection group and $N\unlhd \widetilde{W}$, then one of the following is true:
        \begin{enumerate}[label=(\roman*)]
            \item $N$ is a $W$-invariant sublattice $L\leqslant Q(R^{\vee})$;
            \item $N$ is an extension by $\{\pm1\}$ of a $W$-invariant sublattice $L$ satisfying $2Q(R^{\vee})\leqslant  L\leqslant Q(R^{\vee})$;
            \item There exists a spherical reflection group $W'$, a subgroup $Z\leqslant Z(W')$, and a graph homomorphism $\psi:\widetilde{W}\to W'$ such that $N$ is the kernel of the homomorphism $\widetilde{W}\to W'/Z$ induced by $\psi$.
        \end{enumerate}
        \label{thm: normal subs of W~}
    \end{theorem}
    The proof of theorem \ref{thm: normal subs of W~} involves explicitly constructing the normal subgroups of $\widetilde{W}$ using the correspondence outlined in Remark \ref{rmk: correspondence with triples}, and then observing that each of these subgroups fits into one of the above categories.
    
    \begin{remark}[The correspondence between $N\unlhd \widetilde{W}_n$ and triples $(K,L,c)$]
    \label{rmk: correspondence with triples}
    The set of normal subgroups of $ \widetilde{W}=Q(R^{\vee})\rtimes W$ is in correspondence with the set of triples $(K,L,\bar{c})$ where $L\leqslant Q(R^{\vee})$ is a $W$-invariant sublattice, $K\unlhd W$ acts trivially on $Q(R^{\vee})/L$,  and $\bar{c}:K\to Q(R^{\vee})/L$ is a 1-cocycle satisfying $\overline{c}^w=\bar{c}$ for all $w\in W$. 
    For ease of notation in the third coordinate of $(K,L,\bar{c})$, instead of writing $\bar{c}$, we will often write $(c(x_1),\dots,c(x_k))$ where $c:K\to Q(R^{\vee})$ is any lift of $\bar{c}$ (note that $c$ need not be a homomorphism), and $x_1,\dots,x_k$ is the normally-generating set for $K$ used by Maxwell in Tables 6 and 7. However, when $\bar{c}$ is trivial, we simply write the triple as $(K,L,0)$. 
    
    A normal subgroup $N\unlhd \widetilde{W}$ induces the triple $(\pi(N),N\cap Q(R^{\vee}),\bar{c})$ where $c(x)x\in N$ for all $x\in \pi(N)$. Conversely, given $(K,L,\bar{c})$ the corresponding $N\unlhd \widetilde{W}$ is 
    \begin{equation}
        \label{equ: normal subgroup N~(K,L,c)}
        N=\{(t+c(x),x):t\in L, x\in K\}.
    \end{equation} 
    \end{remark}
    
    Maxwell scales the roots of $W$ so that $\lb \widetilde{\alpha},\widetilde{\alpha}\rb=2$. Hence the lattices Maxwell associates with $W$ differ slightly from those of Bourbaki in {\rm\cite{bourbaki_groupes_1981}}, although this only occurs when $W=C_n$ or $G_2$. Naturally, this difference in the lattices is just in the labelling, not in the structure, but for the sake of clarity we will follow Bourbaki. Hence the labels of normal subgroups of $\widetilde{C}_n$ and $\widetilde{G}_2$ we list differ from those used by Maxwell.

    \begin{remark}[Isomorphism types of quotients of $\widetilde{W}$]
    \label{rmk: isom type small quot}
        Suppose $N=(K,L,(c(x_1),\dots,c(x_m)))$ is a normal subgroup of $\widetilde{W}$ with $|\widetilde{W}:N|\leq |W|$. Then $N$ and quotient group $\widetilde{W}/N$ fit into the following short exact sequences:
        \begin{gather}
            0\to L \to N \to K\to 1 ,\\ 
             0\to Q(R^{\vee})/L\to \widetilde{W}/N \to W/K\to 1. 
             \label{equ: short exact seqs for H norm in W~}
        \end{gather}
        Moreover, these sequences split, so that $N\cong L\rtimes K$ and $\widetilde{W}/N\cong (Q(R^{\vee})/L)\rtimes (W/K)$, if $c\equiv 0$. In general, we use Theorem \ref{thm: normal subs of W~} to understand the quotient of $\widetilde{W}/N$ where $N$ corresponds to the triple $(K,L,c)$. 
        
        The subgroup $N$ is in case $(i)$ of Theorem \ref{thm: normal subs of W~} if and only if  $K=\{1\}$. This implies that $c:K\to Q(R^{\vee})$ is trivial so, in this case, $\widetilde{W}/N$ is isomorphic to $(Q(R^{\vee})/L)\rtimes W$.
        
        Similarly, $N$ is in case $(ii)$ if and only if $K=\{\pm1\}$. 
        We are interested in $N$ with index at most $|W|$ in $\widetilde{W}$ and
        $$|\widetilde{W}:N|=|Q(R^{\vee}):L|\ |W:\{\pm1\}|=
        |Q(R^{\vee}):L|\ \frac{|W|}{2}
        \leq |W|$$ 
        if and only if $|Q(R^{\vee}):L|\leq 2 $; so, in particular, $\Aut(Q(R^{\vee})/L)$ is trivial. 
        If $c$ is trivial, then the short exact sequence in Equation \eqref{equ: short exact seqs for H norm in W~} splits, and since $\Aut(Q(R^{\vee})/L)$ is trivial, we find that
        $$\widetilde{W}/N\cong( Q(R^{\vee})/L ) \rtimes (W/\{\pm1\}) \cong (Q(R^{\vee})/L) \times  W/\{\pm1\} .$$
        If $c(-1)\neq 0$, however, finding the isomorphism type of $\widetilde{W}/N$ is less straightforward. Fortunately, there is only one affine group $\widetilde{W}$ where both $\{\pm1\}\leqslant W$, and which contains an index 2 sublattice $L{\leqslant} Q(R^{\vee})$. This is $\widetilde{W}=\widetilde{C}_n$ with $L=Q(R)$ and, as shown in Table $6$ of {\rm\cite{maxwell_normal_subs}}, there does exist a nontrivial map $\bar{c}$ induced by $c(-1)=\al_n^{\vee}$. In part $(a)$ of Lemma \ref{lem: exceptional quotients of Cn~}, we will show that the quotient of $\widetilde{C}_n$ by $(\{\pm1\},Q(R),\al_n^{\vee})$ is isomorphic to $C_n$.  
        
        Finally, if $N$ is in case $(iii)$, then in Table $7$ Maxwell exhibits $W'$, $\psi:\widetilde{W}\to W'$ and $Z\leqslant Z(W')$, so that $N$ is the kernel of  $\widetilde{W}\to W'/Z$ induced by $\psi$. Hence,
        $$ \widetilde{W}/N\cong  \im(\psi)= W'/Z. $$
        Suppose $W'\cong W_1\times\dots \times W_k$. If $Z\leqslant Z(W_i)$ for some $i$, then without loss of generality $i=1$ and $\widetilde{W}/N$ is isomorphic to $W_1/Z\times W_2\times\dots\times W_k$. In the cases where $Z\nleqslant Z(W_i)$ for every $i$, then the isomorphism type of $\widetilde{W}/Z$ is either computed in $(b)$ of Lemma \ref{lem: exceptional quotients of Cn~}, or is resolved by the observation that $\FC_2^3/\FC_2\cong \FC_2^2$.
    \end{remark}

    \begin{lemma}
    \label{lem: exceptional quotients of Cn~}
        \begin{enumerate}[label=(\alph*)]
             \item Let $n\geq3$, then
             $$ \widetilde{C}_n/(\{\pm1\},Q(R),\al_n^{\vee})\cong C_n.$$
            \item Suppose either that $W_n$ has type $C$ and $n\geq3$, or that $W_n$ has type $D$ and $n\geq3$ is even, then:
            $$ (\FC_2\times W_n)/\lb (-1,-1)\rb\cong W_n. $$
        \end{enumerate}
    \end{lemma}
    \begin{proof}            
        For part $(a)$, define the epimorphism $\phi:\widetilde{C}_n\to C_n$ on $Q(R^{\vee})\rtimes C_n\cong \widetilde{C}_n$ by 
        $$ \phi(v,w)=\chi(v)w ,$$
        where $\chi:Q(R^{\vee})\to \{\pm1\}$ is the unique surjective homomorphism with kernel $Q(R)$. Since $\al_n^{\vee}\notin Q(R), $ we find that
        $$ \ker (\phi)=\{(v,1):v\in Q(R)\}\cup \{(v,-1):v\notin Q(R)\}=(\{\pm1\},Q(R),\al_n^{\vee}).$$

        For part $(b)$, note that $\lb (-1,-1)\rb$ is the kernel of the epimorphism $\FC_2\times W_n\to W_n$ defined by $ ((-1)^{\eps},w)\mapsto z_0^{\eps}w $.
    \end{proof}

\section{The homomorphisms tables of $\Hom(\widetilde{W}_n,W_n)$ for types $A$, $B$, $C$ and $D$}
\label{sec: hom tables ABCD}  

    The symmetry group of the $n$-cube acts on $\RR^n$ by permuting and inverting an orthonormal basis $e_1,\dots,e_n$ of $\RR^n$, and is isomorphic to the irreducible spherical reflection group of type $B_n\cong C_n$. 
    Let $r_{+1,i,j}$ be the reflection with root $e_i-e_j$ which permutes the $i$th and $j$th coordinates, $r_{-1,i,j}$ the reflection with root $e_i+e_j$ which permutes \textit{and} inverts the $i$th and $j$th coordinates, and $r_{i}$ the reflection with root $e_i$ which inverts the $i$th coordinate.    

    The subgroup of coordinate inversions, $(C_n)_2:=\lb r_1,\dots,r_n\rb \cong \FC_2^n$, is normal in $C_n$, and 
    $$C_n=(C_n)_2\rtimes S_+\cong \FC_2\wr \FS_n $$ 
    where $S_+:=\lb r_{+1,1,2},\dots, r_{+1,n-1,n}\rb \cong \mathfrak{S}_n$ is the subgroup of coordinate permutations.
    Since any $w\in C_n$ can be uniquely written as $w=x\sig$ where $x\in (C_n)_2$ and $\sig\in S_+$, we can define the \emph{sign} of $w$ to be the tuple $(\det(x),\det(\sigma))$. This encodes the parity of the number of reflections used to write $x$ and $\sig$.

    We follow the notation for normal subgroups of Coxeter groups, which relates a normal subgroup to the conjugacy class(es) of its generators, used by Maxwell in {\rm\cite{maxwell_normal_subs}}. In particular, we denote by $(C_n)_1$ be the index $2$ subgroup of $C_n$ consisting of isometries which involve an even number of coordinate-inversions. For $n\geq4$, this subgroup is isomorphic to the irreducible spherical reflection group of type $D_n$, and hence isomorphic to $\FC_2^{n-1}\rtimes \FS_n$. The subgroup $S_+\leq C_n$ of coordinate permutations is isomorphic to the irreducible spherical reflection group of type $A_{n-1}$. 

    \begin{theorem}{\rm\cite[p. 21,28]{AutW_Franzsen}}
        For $n\neq 4$, the automorphism group of $C_n$ has order $2^{n+1}n!=2|C_n|$ when $n$ is even; and $2^nn!=|C_n|$ when $n$ is odd. Similarly for $n\neq 4$, the automorphism group of $D_n$ has order $2^{n}n!=2|D_n|$ when $n$ is even; and $2^{n-1}n!=|D_n|$ when $n$ is odd. For $n\neq5$, the automorphism group of $A_n$ has order $(n+1)!=|A_n|$. 
        \label{theorem: |Aut| for Cn, Dn, An}
    \end{theorem}

    Let $w=(x,\sigma)\in C_n$ and suppose $\sigma\in \mathfrak{S}_n$ has cycle-type $(\la_1,\dots,\la_l)$. We can decompose $w$ as $w_1\dots w_l$ so that $w_i$ induces a $\la_i$-cycle in $\mathfrak{S}_n$. We can also partition the orthonormal basis of $\RR^n$ into sets $P_1,\dots,P_l$ with $|P_i|=\la_i$ in such a way that the action of $w$ preserves this partition, and moreover each $P_i$ is pointwise fixed by the action of $w_j$ when $j\neq i$. 
    Then, for any $e_j\in P_i$, we have $w^{k_i}(e_j)=w_i^{k_i}(e_j)=\pm e_j$, and this sign is independent of the choice of $e_j$. We then say that $w_i$ has signed cycle-type $\pm \la_i$, and $w$ has \emph{signed cycle-type} given by the multiset $\{\pm \la_1,\dots,\pm \la_l\}$. For $w\in S_+\cong \FS_n$, the signed cycle type of $w$ coincides with the usual permutation cycle-type, and so determines the conjugacy class of $w$. 

    \begin{theorem}[{\rm{\cite[p. 306]{carter_conjugacy_1972}}}]
        Two elements of $C_n$ are conjugate in $C_n$ if and only if they have the same signed cycle-type. Two elements in $D_n$ are conjugate in $D_n$ if and only if they have the same signed cycle-type, except when every signed cycle has positive sign and even length, in which case there are two conjugacy classes.
    \end{theorem}
    
    If $w$ is an involution, its signed cycle-type consists of $+2$, $-1$, and  $+1$, with some multiplicities. Hence we can assign to $w$ an \emph{involution-type} $(t,u)$, where $t$ encodes the multiplicity of $-1$ in the signed cycle-type of $w$, and $u$ the multiplicity of $+2$. Note that $2t+u\leq n$, and the sign of an involution of involution-type $(t,u)$ is $((-1)^t,(-1)^u)$.

    \begin{example}
        The reflection $r_{\pm1,i,j}$ has signed cycle-type $\{+2,+1,\dots,+1\}$, involution-type $(0,1)$, and sign $(1,-1)$.
        The reflection $r_i$ has signed cycle-type $\{-1,+1,\dots,+1\}$, involution-type $(1,0)$ and sign $(-1,1)$.
        The central element $z_0=r_1\dots r_n$ has signed cycle-type $\{-1,\dots,-1\}$, involution-type $(n,0)$ and sign $((-1)^n,1)$.
        The element $r_1r_{+1,1,2}\in C_2$ has signed cycle-type $\{-2\}$ and sign $(-1,-1)$. 
    \end{example}

\subsection{Images of homomorphisms}
\label{subsec: images of homs ABCD}

    This section counts the number of possible images of homomorphisms $\widetilde{W}_n\to W_n$ when $W_n=A_n$, $B_n$, $C_n$, or $D_n$. Since $B_n\cong C_n$, we will state the Lemmas in this section in terms of $C_n$, but they also apply to $B_n$. We begin by computing $\V_1(W_n)$ which we use to count the number of subgroups of $W_n$ isomorphic to $\FC_2$.

    \begin{lemma}{\rm\cite[Lemma 3.5]{random_endos}}
        Suppose $W_n$ is $A_{n-1}$, $C_n$, or $D_n$ where $n\geq2$. The number of elements of $W_n$ with order at most $2$ is
        $$ \mathcal{V}_1(W_n)=\frac{2^{\al(n)}n!}{\lfloor \frac{n}{2}\rfloor!}+\sum_{k=0}^{\lfloor \frac{n}{2}\rfloor-1} \frac{|W_n|}{2^{\tau k}\;k!\;(n-2k)!} $$
        where $\tau=2$ and $\al(n)=\frac{1-(-1)^n}{2}$ when $W_n=C_n$, $\tau=2$ and $\al(n)=0$ when $W_n=D_n$, and $\tau=1$ and $\al(n)=-\lfloor \frac{n}{2}\rfloor$ when $W_n=A_{n-1}$. 
        \label{lem: V1 for Cn, Dn, An} 
    \end{lemma}
    
    Lemma \ref{lem: commuting pairs and triples in Cn} will be used to find the number of subgroups of $C_n$ isomorphic to $\FC_2^2$ and $\FC_2^3$. The proof of this Lemma proceeds by first observing that a $k$-tuple of pairwise commuting elements with order at most $2$ induces a partition on $\{1,\dots,n\}$ into sets of size $1,2,\dots,2^k$, before counting the number of tuples which induce said partition, and the number of such partitions.

    \begin{lemma}
        When $n\geq2$, 
        the number of pairs of commuting elements in $C_n$ with order at most $2$ is 
        $$\V_2(C_n)=\sum_{k=0}^{\lfloor \frac{n}{4}\rfloor}\; \sum_{l=0}^{\lfloor \frac{n}{2}\rfloor-2k} \frac{(\frac{3}{8})^l2^{2n}n!}{2^{7k} \ k!\ l!\ (n-4k-2l)!} . $$
        Furthermore, the number of ordered triples of pairwise commuting elements in $C_n$ with order at most $2$ is 
        \begin{equation*}
        \V_3(C_n)=\sum_{k=0}^{\lfloor \frac{n}{8}\rfloor} \;
        \sum_{l=0}^{\lfloor \frac{n}{4}\rfloor-2k} \;
        \sum_{m=0}^{\lfloor \frac{n}{2}\rfloor-4k-2l} 
        \frac{2^{3n}7^{l+m}n!}{2^{20k+10l+4m} \ k!\ l!\ m!\ (n-8k-4l-2m)!}    .   
        \end{equation*}
    \label{lem: commuting pairs and triples in Cn}
    \end{lemma}
    \begin{proof}
        For $\V_2(C_n)$, see Lemma 3.6 of {\rm\cite{random_endos}}.
        
        Fix an element $x$ with order at most $2$ and involution-type $(t,u)$ in $C_n$. Then $x$ has the form 
        $$ x=r_{a_1}\dots r_{a_t}r_ {\eps_1,b_1,c_1}\dots r_{\eps_u,b_u,c_u}  $$
        where $\eps_i\in\{\pm1\}$ and $a_i,b_i,c_i\in \{1,\dots,n\}$ are all distinct. Let $\mathrm{Fix}(x)$ be the set of indices of the orthonormal basis vectors of $\RR^n$ which are fixed by the action of $x$ on $\RR^n$. Suppose $y$ and $z$ are commuting elements of $C_n$ with order at most $2$ which commute with $x$. Then $y$ and $z$ are products of terms with the following forms:
        \begin{enumerate}[label=(\roman*)]
            \item $r_{\eps_l,b_i,c_i}$ where $\eps_l\in\{\pm1\}$,
            \item $r_{\eps_l,a_i,a_j}$ where $\eps_l\in\{\pm1\}$,
            \item $r_{\eps_l,d,d'}$ where $d,d'\in\mathrm{Fix}(x)$ and $\eps_l\in\{\pm1\}$,
            \item $r_{\eps_k,b_i,b_j}r_{\eps_l,c_i,c_j}$ where $\eps_k,\eps_l\in\{\pm1\}$ so that $\eps_i\eps_j=\eps_k\eps_l$,
            \item $r_{b_i}r_{c_i}$,
            \item $r_{a_i}$,
            \item $r_d$ where $d\in\mathrm{Fix}(x) $.
        \end{enumerate}
        Given $x$, $y$ and $z$, one can orthogonally decompose $\RR^n$ into subspaces of dimension $1$, $2$, $4$, or $8$ which are preserved under the action of $\lb x,y,z\rb$ on $\RR^n$, and for which the dimensions of these subspaces are minimal with respect to this property. In particular, these subspaces are spanned by disjoint subsets of the orthonormal basis $\{e_1,\dots,e_n\}$ of $\RR^n$, and so this decomposition induces a partition of $\{1,\dots,n\}$ into parts of size $1$, $2$, $4$, or $8$. Given such a decomposition of $\RR^n$, we will count how many triples $(x,y,z)$ of pairwise commuting elements with order at most $2$ produce it. Suppose $P=\{P_{\al}\}$ is a partition of $\{1,\dots,n\}$ with parts of size $1$, $2$, $4$, and $8$.  If $x$, $y$ and $z$ are pairwise commuting elements of $C_n$ with order at most $2$ which induce $P$, then they are products of terms of the forms listed in Table \ref{tab: terms in commuting refls}. 
        
        \begin{table}[h] 
        \begin{adjustbox}{width=\textwidth}
        \begin{tabular}{c | c c c | c c c}
                $P_{\alpha}$  & Term in $x$ & Term in $y$ & Term in $z$ & Allocation & Coefficient & Permutation  \\
             & & & & choices & choices & choices\\
            \hline 
            \hline
                $\{i\}$ & $r_i^{\ka_1}$ & $r_i^{\ka_2}$ & $r_i^{\ka_3}$ & 1 & $2^3$  & 1 \\
            \hline 
                & $r_{\eps_1,i,j} $&$r_{\eps_2,i,j}$ & $r_{\eps_3,i,j}$  &1 & $2^3$  & 1 \\
                $\{i,j\}$ & $r_{\eps_1,i,j} $& $r_{\eps_2,i,j}$ &$(r_ir_j)^{\ka}$ & 1 & $2^3$ & 3 \\
                & $r_{\eps,i,j}$ & $(r_ir_j)^{\ka_1} $&$(r_ir_j)^{\ka_2}$ & 1 & $2^3$ & 3 \\
            \hline 
                & $r_{\eps_1,i,j}r_{\eps_2,k,l}$ & $r_{\eps_3,i,k}r_{\eps_1\eps_2\eps_3,j,l}$ & $r_{\eps_4,i,l}r_{\eps_1\eps_2\eps_4,j,k}$ & $2\cdot 3$ & $2^4$ & 1 \\ 
                $\{i,j,k,l\}$ & $r_{\eps_1,i,j}r_{\eps_2,k,l}$ & $r_{\eps_3,i,k}r_{\eps_1\eps_2\eps_3,j,l}$ & $r_{\eps_4,i,k}r_{\eps_1\eps_2\eps_4,j,l}$ & $2\cdot 3$ & $2^4$ & 3 \\ 
                & $r_{\eps_1,i,j}r_{\eps_2,k,l}$ & $r_{\eps_3,i,k}r_{\eps_1\eps_2\eps_3,j,l}$ & $(r_ir_jr_kr_l)^{\ka}$ & $2\cdot 3$ & $2^4$ & 3 \\ 
            \hline 
                \multicolumn{1}{l|}{$\{i,j,k,l,\ldots$} & 
                \multicolumn{1}{l}{$r_{\eps_1,i,j}r_{\eps_2,k,l}\ldots$} &
                \multicolumn{1}{l}{$r_{\eps_5,i,k}r_{\eps_1\eps_2\eps_5,j,l}\ldots$ }& 
                \multicolumn{1}{l|}{$r_{\eps_7,i,m}r_{\eps_1\eps_2\eps_7,k,n}\ldots$ }& $2^4\cdot3^2\cdot5\cdot7$ & $2^7$ & 1\\
                \multicolumn{1}{r|}{$m,n,o,p\}$} & $r_{\eps_3,m,n}r_{\eps_4,o,p}$ & $r_{\eps_6,m,o}r_{\eps_3\eps_4\eps_6,n,p}$ & $r_{\eps_3\eps_4\eps_6,l,o}r_{\eps_2\eps_4\eps_5\eps_6\eps_7,m,p}$ & & & \\
        \end{tabular}
        \end{adjustbox}
        \caption{Terms that appear in triples of pairwise commuting elements of $C_n$ with order at most 2 with respect to the partite sets they induce. The final three columns count the choices for each triple of terms. Given the set $P_{\al}$, the fifth column counts how its elements can be distributed, the sixth column counts the choices for $\eps\in \{\pm1\}$ and $\ka\in \{0,1\}$, and the last column accounts for asymmetry between $x$, $y$ and $z$.}
        \label{tab: terms in commuting refls}
        \end{table}
        
        From Table \ref{tab: terms in commuting refls}, we see that for a given part $P_{\al}$, if $|P_{\al}|=1$ (respectively $|P_{\al}|=2$, $4$ or $8$), there are $2^3$ (resp. $2^3\cdot 7=2^2\cdot 7\cdot2!$, $2^5\cdot 3\cdot 7=2^2\cdot 7\cdot 4!$, or $2^{11}\cdot 3^2\cdot 5 \cdot 7=2^4\cdot 8!$) ways that $P_{\al}$ can index terms in $x$ and $y$.
    
        There are 
        $\frac{N!}{(L!)^K\ K!\ (N-KL)!}$
        ways to choose $K$ order-$L$ subsets from $N$ elements. A partition $P$ has $0\leq k \leq \lfloor \frac{n}{8}\rfloor $ parts $P_{\al}$ of order $8$ from $n$ elements, $0\leq l\leq \lfloor\frac{n-8k}{4}\rfloor$ subsets of order $4$ from the remaining $n-8k$ elements, $0\leq m \leq \lfloor\frac{n-8k-4l}{2}\rfloor$ parts of order $2$ from $n-8k-4l$ elements, and there are $n-8k-4l-2m$ elements in singleton parts.
        Hence the total number of ordered triples of pairwise commuting elements with order at most 2 is:
        \begingroup\fontsize{10}{11}\selectfont 
        \begin{align*}
            \sum_{k=0}^{\lfloor \frac{n}{8}\rfloor} 
            \frac{(2^4\cdot 8!)^k\ n!}{(8!)^k\ k!\ (n-8k)!}
            \sum_{l=0}^{\lfloor \frac{n-8k}{4}\rfloor}  
            \frac{(2^2\cdot 7\cdot 4!)^l\ (n-8k)!}{(4!)^l\ l!\ (n-8k-4l)!} 
            \sum_{m=0}^{\lfloor \frac{n-8k-4l}{2}\rfloor} 
            \frac{(2^2\cdot 7\cdot2!)^m(n-8k-4l)!}{(2!)^m\ m!\ (n-8k-4l-2m)!}\
            2^{3(n-8k-4l-2m)}   .
        \end{align*}
        \endgroup
    \end{proof}

    We have computed the number of ordered pairs and triples of commuting elements with order at most 2, and Lemma \ref{lem: subs of Cn isom (Z/2Z)^2 and (Z/2Z)^3} uses this to find the number of subgroups of $C_n$ isomorphic to $\FC_2^2$ and $\FC_2^3$. 
    
    \begin{lemma}
        Let $n\geq2$. 
        \begin{enumerate}[label=(\alph*)]
            \item The number of subgroups of $C_n$ isomorphic to $\FC_2^2$ is 
        $$\frac{1}{6}(\V_2(C_n)-3\V_1(C_n)+2)
        .$$
            \item The number of subgroups of $C_n$ isomorphic to $\FC_2^3$ is
        $$  \frac{1}{168} (\V_3(C_n)-7\V_2(C_n)+14\V_1(C_n)-8).$$
        \end{enumerate}         
    \label{lem: subs of Cn isom (Z/2Z)^2 and (Z/2Z)^3}
    \end{lemma}
    \begin{proof}
        Part $(a)$ is Lemma 3.7 of {\rm\cite{random_endos}}. For part $(b)$, note that subgroups of $C_n$ isomorphic to $\FC_2^3$ are generated by triples of distinct, pairwise commuting involutions. We will remove from $\V_3(C_n)$ the number of triples which do not generate a subgroup isomorphic to $\FC_2^3$, but instead generate subgroups isomorphic to $\FC_2^2$, $\FC_2$ or $\FC_1$. 
        For each subgroup $H\leqslant C_n$ isomorphic to $\FC_2^2$ (respectively $\FC_2$ or $\FC_1$), there are $42$ (resp. $7$ or $1$) triples $(x,y,z)\in H^3$ which generate $H$, and each of these is a triple of pairwise commuting elements with order at most $2$. Hence, by part $(a)$ of Lemma \ref{lem: commuting pairs and triples in Cn} and the fact that $|\Aut(\FC_2^3)|=168$, the number of subgroups isomorphic to $\FC_2^3$ is 
        $$ \frac{1}{168}\big(\V_3(C_n)-\frac{42}{6}(\V_2(C_n)-3\V_1(C_n)+2)-7(\V_1(C_n)-1)-1\big)=\frac{1}{168} \big(\V_3(C_n)-7\V_2(C_n)+14\V_1(C_n)-8\big). $$
    \end{proof}    
    
    We will count the number of subgroups of $C_n$ which are isomorphic to $\FC_2^k\times \mathfrak{S}_n$ for $k\geq0$. Note that $S_+:=\lb r_{+1,1,2},\dots, r_{+1,n-1,n}\rb $ and $S_-:=\lb z_0r_{+1,1,2},\dots, z_0r_{+1,n-1,n}\rb $ are both subgroups of $C_n$ isomorphic to $\FS_n$.

    \begin{lemma}
        Suppose $n\geq 3$ and $n\ne 4$. 
        \begin{enumerate}[label=(\alph*)]
            \item There are $2^n$ subgroups of $C_n$ isomorphic to $\FS_n$, half of which are conjugate in $C_n$ to $S_+$, and the other half to $S_-$. 
            \item When $n$ is odd, there are $2^{n-1}$ subgroups of $D_n$ isomorphic to $\FS_n$, and these are all conjugate in $D_n$ to $S_+$. When $n$ is even, there are $2^{n}$  subgroups of $D_n$ isomorphic to $\FS_n$, and $\Aut(D_n)$ acts transitively on the this set of subgroups.
            \item There are $2^{n-1}$ subgroups of $C_n$ isomorphic to $\FC_2\times\FS_n$, all of which are conjugate in $C_n$ to $\lb z_0, S_+\rb$.
            \item There is no subgroup of $C_n$ isomorphic to $\FC_2^k\times \FS_n$ for $k\geq 2$.
        \end{enumerate}
        \label{lem: medium size images BCD}
    \end{lemma}
    \begin{proof}
        For parts $(a)$, $(b)$ and $(c)$ see Lemmas 3.8 and 3.9 of {\rm\cite{random_endos}}. Part $(d)$ follows from $(a)$, $(c)$, and the fact that $z_0$ is central. Together these imply that $z_0$ is the unique involution in $C_n$ which commutes with any subgroup isomorphic to $\FS_n$.
    \end{proof}

    Finally, in Lemma \ref{lem: large images BCD} we examine the subgroups $H\leqslant W_n$ for which $|W_n:H|\leq 2$. Since these subgroups are necessarily normal and there are finitely many normal subgroups of $W_n$, this reduces to checking the isomorphism types of these subgroups.
    
    \begin{lemma}
    \label{lem: large images BCD}
        Suppose $n\geq 3$ and $n\ne 4$.
        \begin{enumerate}[label=(\alph*)]
            \item If $n$ is odd, there are two subgroups, $(C_n)_1$ and $C_n^+$, of $C_n$ isomorphic to $C_n/\{\pm1\}$. Further, the automorphism group of $C_n/\{\pm1\}$ has order $2^{n-1}n!$. When $n$ is even, there is no subgroup of $C_n$ isomorphic to $C_n/\{\pm1\}$. 
            \item Suppose $n$ is even, there is no subgroup of $C_n$ or $D_n$ isomorphic to either $D_n/\{\pm1\}$ or $\FC_2\times D_n/\{\pm1\}$.
            \item When $n$ is odd, there are two subgroups, $(C_n)_1$ and $C_n^+$, of $C_n$ isomorphic to $D_n$. When $n$ is even, $(C_n)_1$ is the only subgroup of $C_n$ isomorphic to $D_n$.
            \item When $n$ is odd there is one subgroup, $\lb z_0\rb\times (C_n)_1=C_n$, of $C_n$ isomorphic to both $\FC_2\times D_n$ and $\FC_2\times C_n/\{\pm1\}$, and its automorphism group has order $2^nn!$. When $n $ is even, there is no subgroup of $C_n$ isomorphic to either $\FC_2\times D_n$ or $\FC_2\times C_n/\{\pm1\}$. 
        \end{enumerate}
    \end{lemma}
    \begin{proof}
        Parts $(a)$ and $(b)$ follow from the proof of Lemma 3.10 of {\rm\cite{random_endos}}. For part $(c)$, Recall that $(C_n)_1$ is isomorphic to $D_n$. Part $(a)$ completes the proof when $n$ is odd. Suppose $n$ is even, it remains to show that $C_n+$ and $\prescript{}{+}{C_n}$ are not isomorphic to $D_n$. 
        The subgroups $C_n^+$ and $\prescript{}{+}{C_n}$ consist of the elements of $C_n$ with sign $\pm(+1,+1)$ and $(\pm1,+1)$, respectively. Hence neither of them contains $S_{\pm}$, since these subgroups contain  either $r_{+1,1,2}$ or $z_0r_{+1,1,2}$, and these elements have sign $(+1,-1)$. Lemma \ref{lem: medium size images BCD} part $(a)$ together with the fact that $C_n^+,\prescript{}{+}{C_n}\unlhd C_n$ imply that neither $C_n^+$ nor $\prescript{}{+}{C_n}$ contain any subgroup isomorphic to $\FS_n$. Then, because $D_n$ has $2^n$ subgroups isomorphic to $\FS_n$, these are not isomorphic to $D_n$.

        With regard to part $(d)$ when $n$ is odd, $z_0\notin (C_n)_1$, so $\lb z_0\rb \times (C_n)_1=C_n$ is isomorphic to both $\FC_2\times D_n$ and $\FC_2\times C_n/\{\pm1\}$ by part $(a)$. When $n$ is even, the claim follows from parts $(a)$ and $(c)$, and the fact that $z_0\in (C_n)_1$ is the only involution which commutes with $S_{+}\leqslant (C_n)_1$.  
    \end{proof}

\subsection{Counting $\Hom(\widetilde{W}, W)$ for types $A$, $B$, $C$ and $D$}
\label{subsec: counting homs ABCD}

    In this section we prove parts $(a)$-$(c)$ of Theorem \ref{thm: main theorem} by counting the number of homomorphisms $\widetilde{W}_n\to W_n$ for $W_n=A_n$, $B_n$, $C_n$, or $D_n$. Suppose $W_n$ is $A_n$, $B_n$, $C_n$ or $D_n$. Theorem \ref{theorem: |Aut| for Cn, Dn, An} states the order of the automorphism group of $W_n$, Section \ref{thm: normal subs of W~} describes the normal subgroups $N\unlhd \widetilde{W}_n$ and the quotient group $\widetilde{W}_n/N$, and Lemmas \ref{lem: V1 for Cn, Dn, An}-\ref{lem: large images BCD} count the number of subgroups of $W_n$ isomorphic to $W_n/N$ and, where applicable, state the order of $\Aut(\widetilde{W}_n/N)$. These data are collated in Table \ref{tab: counting Hom(W~,W) for A,B,C,D}.

    \begingroup\fontsize{8}{10}\selectfont 
    \begin{longtable}{ m{0.4in} m{1.3in} m{0.72in} m{1.3in} m{0.85in} m{1.1in}  }
        $ |\widetilde{W}:N|$ & $(K,L,(c(x_1),\dots,c(x_k)))$\newline \mbox{}\hfill$\sim N\unlhd \widetilde{W} $\hfill\mbox{} &  $\widetilde{W}/N$ &$|\{X\leqslant W:X\cong \widetilde{W}/N\}|$ & $|\Aut(\widetilde{W}/N)|$ & $|\{\varphi\in \mathrm{Hom}(\widetilde{W},W):$\newline \mbox{}\hfill$\ker\varphi=N\}|$ \\ 
    \noalign{\vskip 4pt}\hline\hline\noalign{\vskip 4pt}
        \multicolumn{6}{l}{$A_n$ $n\ge 4$ and $n\neq5$:}\\ 
            $1$ & $(A_n,Q(R),0)$& $\{1\}$ & $1$ & $1$ & $1$ \\ 
             
            $2$ & $(A_n^+,Q(R),0)$ & $\FC_2$ & $\V_1(A_n)-1$  & $1$ & $\V_1(A_n)-1$ \\   

            $(n+1)!$ & $(\{1\},Q(R),0)$  & $A_n$ & $1$ & $(n+1)!$ & $(n+1)!$ \\
    \noalign{\vskip 4pt}\hline\noalign{\vskip 4pt}
        \multicolumn{6}{l}{$B_n$ $n\ge 3$ and $n\neq 4,6$:}\\ 
            $1$ & $(B_n,Q(R^{\vee}),0)$& $\{1\}$ & $1$ & $1$ & $1$ \\ 
             
            $2$ & $(B_n^+,Q(R^{\vee}),0)$ & $\FC_2$ & $\V_1(B_n)-1$  & $1$ & $\V_1(B_n)-1$ \\             
             
            $2$ & $((B_n)_1,Q(R^{\vee}),0) $ &  $\FC_2 $ & $\V_1(B_n)-1$ & $1$ & $\V_1(B_n)-1$ \\
             
            $2$ & $(\prescript{}{+}{B_n},Q(R^{\vee}),0) $ &  $\FC_2 $  & $\V_1(B_n)-1$ & $1$ & $\V_1(B_n)-1$ \\
             
            $4$ & $((B_n)_1^+ ,Q(R^{\vee}),0)$ & $\FC_2^2 $   & $\frac{1}{6}(\V_2(B_n)-3\V_1(B_n))+2)$& $6$ & $\V_2(B_n)-3\V_1(B_n)+2$ \\
             
            $n!$ & $((B_n)_2,Q(R^{\vee}),0)$  & $\mathfrak{S}_n$ & $2^n $& $n!$ & $2^nn!$\\
             
            $2 \cdot n!$& $((B_n)_2^+,Q(R^{\vee}),0)$ & $\FC_2\times \mathfrak{S}_n$  & $2^{n-1}$& $2\cdot n! $& $2^nn!$ \\
                         
            $ 2^{n-1}n! $& $(\{\pm1\},Q(R^{\vee}),0)$ & $B_n/\{\pm1\}$ & 
            $\begin{cases}
                2 & \text{$n$ odd,}\\ 
                0 & \text{$n$ even}.
            \end{cases}$ &
            $\begin{cases}
                2^{n-1}n! & \text{$n$ odd,}\\
                \mathrm{N/A} & \text{$n$ even.}
            \end{cases}$
            & $\begin{cases}
                2^{n}n! & \text{$n$ odd,}\\
                0 & \text{$n$ even.}
            \end{cases}$\\[15pt]

            $ 2^{n-1}n! $& $((B_n)_2,2Q(R),0)$ & $D_n$ & 
            $\begin{cases}
                2 & \text{$n$ odd,}\\ 
                1 & \text{$n$ even}.
            \end{cases}$ &
            $\begin{cases}
                2^{n-1}n! & \text{$n$ odd,}\\
                2^nn! & \text{$n$ even.}
            \end{cases}$
            & $2^nn!$\\[15pt]

            $2^nn!$ & $(\{1\},Q(R^{\vee}),0)$  &  $B_n$ & $1$ & $\begin{cases}
                2^nn! & \text{$n$ odd,} \\ 
                2^{n+1}n! & \text{$n$ even.} 
            \end{cases}$ & 
            $\begin{cases}
                2^nn! & \text{$n $ odd,}\\
                2^{n+1}n! & \text{$n$ even.}
            \end{cases}$\\[15pt]

            $2^nn!$ & $((B_n)_2^+,2Q(R),0)$  &  $\FC_2\times D_n$ & 
            $\begin{cases}
                1 & \text{$n$ odd,} \\ 
                0 & \text{$n$ even.}           
            \end{cases}$ & 
            $\begin{cases}
                2^nn! & \text{$n$ odd,} \\ 
                \mathrm{N/A} & \text{$n$ even.} 
            \end{cases}$ & 
            $\begin{cases}
                2^nn! & \text{$n $ odd,}\\
                0 & \text{$n$ even.}
            \end{cases}$\\[15pt]

            $2^nn!$ & $((B_n)_2^+,2Q(R),\al_{n-1})$  &  $\ B_n$ & $1$ & 
            $\begin{cases}
                2^nn! & \text{$n$ odd,} \\ 
                2^{n+1}n! & \text{$n$ even.} 
            \end{cases}$ & 
            $\begin{cases}
                2^nn! & \text{$n $ odd,}\\
                2^{n+1}n! & \text{$n$ even.}
            \end{cases}$\\[10pt]
    
        \multicolumn{6}{l}{$B_n$ ($n$ odd):}\\ 
            $2^{n-1}n!$  & $((B_n)_2,2Q(R),2\omega_n+\al_n)$  & $B_n/\{\pm1\}$ & $2 $& $2^{n-1}n!$ & $2^nn!$\\[10pt]

        \multicolumn{6}{l}{$B_n$ ($n$ even):}\\ 
            $2^{n-2}n!$  & $((B_n)_2,2P(R),0)$  & $D_n/\{\pm1\}$ & $0 $& $\mathrm{N/A}$ & $0$\\

            $2^{n-1}n!$  & $((B_n)_2,2Q(R),2\omega_n)$  & $D_n$ & $1 $& $2^nn!$ & $2^nn!$\\
             
            $2^{n-1}n!$  & $((B_n)_2^+,2P(R),0)$  & $\FC_2\times D_n/\{\pm1\}$ & $0 $& N/A & $0$\\

            $2^{n-1}n!$  & $((B_n)_2^+,2P(R),\al_{n-1})$  & $B_n/\{\pm1\}$ & $0$ & $\mathrm{N/A}$ & $0$\\
    \noalign{\vskip 4pt}\hline\noalign{\vskip 4pt}
        \multicolumn{6}{l}{$C_n$ $n\ge3$ and $n\neq 4,6$:}\\ 
            $1$ & $(C_n,Q(R^{\vee}),0)$& $\{1\}$ & $1$ & $1$ & $1$ \\ 
             
            $2$ & $(C_n^+,Q(R^{\vee}),0)$ & $\FC_2$ & $\V_1(C_n)-1$  & $1$ & $\V_1(C_n)-1$ \\             
             
            $2$ & $((C_n)_1,Q(R^{\vee}),0) $ &  $\FC_2 $ & $\V_1(C_n)-1$ & $1$ & $\V_1(C_n)-1$ \\
             
            $2$ & $(\prescript{}{+}{C_n},Q(R^{\vee}),0) $ &  $\FC_2 $  & $\V_1(C_n)-1$ & $1$ & $\V_1(C_n)-1$ \\

            $2$ & $(C_n,Q(R),(0,0))$& $\FC_2$ & $\V_1(C_n)-1$ & $1$ & $\V_1(C_n)-1$ \\ 

            $2$ & $(C_n,Q(R),(0,\al_n^{\vee}))$& $\FC_2$ & $\V_1(C_n)-1$ & $1$ & $\V_1(C_n)-1$ \\ 

            $2$ & $(C_n,Q(R),(\al_n^{\vee},0))$& $\FC_2$ & $\V_1(C_n)-1$ & $1$ & $\V_1(C_n)-1$ \\ 

            $2$ & $(C_n,Q(R),(\al_n^{\vee},\al_n^{\vee}))$& $\FC_2$ & $\V_1(C_n)-1$ & $1$ & $\V_1(C_n)-1$ \\ 
             
            $4$ & $((C_n)_1^+ ,Q(R^{\vee}),0)$ & $\FC_2^2 $   & $\frac{1}{6}(\V_2(C_n)-3\V_1(C_n)+2)$ & $6$ & $\V_2(C_n)-3\V_1(C_n)+2$ \\

            $4$ & $((C_n)_1 ,Q(R),0)$ & $\FC_2^2 $   & $\frac{1}{6}(\V_2(C_n)-3\V_1(C_n)+2)$ & $6$ & $\V_2(C_n)-3\V_1(C_n)+2$ \\

            $4$ & $((C_n)_1 ,Q(R),\al_n^{\vee})$ & $\FC_2^2 $   & $\frac{1}{6}(\V_2(C_n)-3\V_1(C_n)+2)$ & $6$ & $\V_2(C_n)-3\V_1(C_n)+2$ \\

            $4$ & $(\prescript{}{+}{C_n},Q(R),(0,0))$ & $\FC_2^2 $   & $\frac{1}{6}(\V_2(C_n)-3\V_1(C_n)+2)$ & $6$ & $\V_2(C_n)-3\V_1(C_n)+2$ \\

            $4$ & $(\prescript{}{+}{C_n},Q(R),(0,\al_n^{\vee}))$ & $\FC_2^2 $   & $\frac{1}{6}(\V_2(C_n)-3\V_1(C_n)+2)$ & $6$ & $\V_2(C_n)-3\V_1(C_n)+2$ \\
            
            $4$ & $(C_n^+ ,Q(R),0)$ & $\FC_2^2 $   & $\frac{1}{6}(\V_2(C_n)-3\V_1(C_n)+2)$& $6$ & $\V_2(C_n)-3\V_1(C_n)+2$ \\

            $4$ & $(C_n^+ ,Q(R),\al_n^{\vee})$ & $\FC_2^2 $   & $\frac{1}{6}(\V_2(C_n)-3\V_1(C_n)+2)$& $6$ & $\V_2(C_n)-3\V_1(C_n)+2$ \\

            $8$ & $((C_n)_1^+ ,Q(R),0)$ & $\FC_2^3 $   & 
            $\frac{1}{168}(\V_3(C_n)-7\V_2(C_n)$
            \newline \mbox{}\hfill$+14\V_1(C_n)-8)$& $168$ 
            & $\V_3(C_n)-7\V_2(C_n)$
            \newline \mbox{}\hfill$+14\V_1(C_n)-8$ \\

            $n!$ & $((C_n)_2,Q(R^{\vee}),0)$  & $\mathfrak{S}_n$ & $2^n $& $n!$ & $2^nn!$\\

            $2 \cdot n!$& $((C_n)_2^+,Q(R^{\vee}),0)$ & $\FC_2\times \mathfrak{S}_n$  & $2^{n-1}$& $2\cdot n! $& $2^nn!$ \\

            $2 \cdot n!$& $((C_n)_2,Q(R),0)$ & $\FC_2\times \mathfrak{S}_n$  & $2^{n-1}$& $2\cdot n! $& $2^nn!$ \\

            $2 \cdot n!$& $((C_n)_2,Q(R),\al_n^{\vee})$ & $\FC_2\times \mathfrak{S}_n$  & $2^{n-1}$& $2\cdot n! $& $2^nn!$ \\

            $4 \cdot n!$& $((C_n)_2^+,Q(R),0)$ & $\FC_2^2\times \mathfrak{S}_n$  & $0$& N/A & $0$ \\
                                    
            $ 2^{n-1}n! $& $(\{\pm1\},Q(R^{\vee}),0)$ & $C_n/\{\pm1\}$ & 
            $\begin{cases}
                2 & \text{$n$ odd,}\\ 
                0 & \text{$n$ even}.
            \end{cases}$ &
            $\begin{cases}
                2^{n-1}n! & \text{$n$ odd,}\\
                \mathrm{N/A} & \text{$n$ even.}
            \end{cases}$
            & $\begin{cases}
                2^{n}n! & \text{$n$ odd,}\\
                0 & \text{$n$ even.}
            \end{cases}$\\[15pt]

            $ 2^{n-1}n! $& $((C_n)_2,2P(R^{\vee}),0)$ & $C_n/\{\pm1\}$ & 
            $\begin{cases}
                2 & \text{$n$ odd,}\\ 
                0 & \text{$n$ even}.
            \end{cases}$ &
            $\begin{cases}
                2^{n-1}n! & \text{$n$ odd,}\\
                \mathrm{N/A} & \text{$n$ even.}
            \end{cases}$
            & $\begin{cases}
                2^{n}n! & \text{$n$ odd,}\\
                0 & \text{$n$ even.}
            \end{cases}$\\[15pt]

            $ 2^{n-1}n! $& $((C_n)_2,2P(R^{\vee}),\alpha_{n})$ & $C_n/\{\pm1\}$ & 
            $\begin{cases}
                2 & \text{$n$ odd,}\\ 
                0 & \text{$n$ even}.
            \end{cases}$ &
            $\begin{cases}
                2^{n-1}n! & \text{$n$ odd,}\\
                \mathrm{N/A} & \text{$n$ even.}
            \end{cases}$
            & $\begin{cases}
                2^{n}n! & \text{$n$ odd,}\\
                0 & \text{$n$ even.}
            \end{cases}$\\[15pt]

            $2^nn!$ & $(\{1\},Q(R^{\vee}),0)$  &  $C_n$ & $1$ & $\begin{cases}
                2^nn! & \text{$n$ odd,} \\ 
                2^{n+1}n! & \text{$n$ even.} 
            \end{cases}$ & 
            $\begin{cases}
                2^nn! & \text{$n $ odd,}\\
                2^{n+1}n! & \text{$n$ even.}
            \end{cases}$\\[15pt]

            $2^nn!$ & $(\{\pm1\},Q(R),0)$ &  $\FC_2\times C_n/\{\pm1\}$ & 
            $\begin{cases}
                1 & \text{$n $ odd,}\\
                0 & \text{$n$ even.}                
            \end{cases}$& 
            $\begin{cases}
                2^nn! & \text{$n$ odd,} \\ 
                \mathrm{N/A} & \text{$n$ even.} 
            \end{cases}$ & 
            $\begin{cases}
                2^nn! & \text{$n $ odd,}\\
                0 & \text{$n$ even.}
            \end{cases}$\\[15pt]

            $2^nn!$ & $(\{\pm1\},Q(R),\al_n^{\vee})$ &  $C_n$ & $1$ & $\begin{cases}
                2^nn! & \text{$n$ odd,} \\ 
                2^{n+1}n! & \text{$n$ even.} 
            \end{cases}$ & 
            $\begin{cases}
                2^nn! & \text{$n $ odd,}\\
                2^{n+1}n! & \text{$n$ even.}
            \end{cases}$\\[15pt]

            $2^nn!$ & $((C_n)_2^+,2P(R^{\vee}),0)$  &  $\FC_2\times C_n/\{\pm1\}$ & 
            $\begin{cases}
                1 & \text{$n$ odd,} \\ 
                0 & \text{$n$ even.}           
            \end{cases}$ & 
            $\begin{cases}
                2^nn! & \text{$n$ odd,} \\ 
                \mathrm{N/A} & \text{$n$ even.} 
            \end{cases}$ & 
            $\begin{cases}
                2^nn! & \text{$n $ odd,}\\
                0 & \text{$n$ even.}
            \end{cases}$\\[15pt]

            $2^nn!$ & $((C_n)_2^+,2P(R^{\vee}),\alpha_{n-1}^{\vee{}})$  &  $\FC_2\times C_n/\{\pm1\}$ &  
            $\begin{cases}
                1 & \text{$n $ odd,}\\
                0 & \text{$n$ even.}                
            \end{cases}$& 
            $\begin{cases}
                2^nn! & \text{$n$ odd,} \\ 
                \mathrm{N/A} & \text{$n$ even.} 
            \end{cases}$ & 
            $\begin{cases}
                2^nn! & \text{$n $ odd,}\\
                0 & \text{$n$ even.}
            \end{cases}$\\[15pt]
            
             $2^nn!$ & $((C_n)_2,2Q(R^{\vee}),\al_n^{\vee})$  &  $\ C_n$ & $1$ &  
            $\begin{cases}
                2^nn! & \text{$n$ odd,} \\ 
                2^{n+1}n! & \text{$n$ even.} 
            \end{cases}$ & 
            $\begin{cases}
                2^nn! & \text{$n $ odd,}\\
                2^{n+1}n! & \text{$n$ even.}
            \end{cases}$\\[15pt]
            
            $2^nn!$ & $((C_n)_2,2Q(R^{\vee}),0)$  &  $\ C_n$ & $1$ &  
            $\begin{cases}
                2^nn! & \text{$n$ odd,} \\ 
                2^{n+1}n! & \text{$n$ even.} 
            \end{cases}$ & 
            $\begin{cases}
                2^nn! & \text{$n $ odd,}\\
                2^{n+1}n! & \text{$n$ even.}
            \end{cases}$\\[15pt]
            
            $2^nn!$ & $((C_n)_2,2Q(R^{\vee}),2\omega_n^{\vee})$  &  $\ C_n$ & $1$ &  
            $\begin{cases}
                2^nn! & \text{$n$ odd,} \\ 
                2^{n+1}n! & \text{$n$ even.} 
            \end{cases}$ & 
            $\begin{cases}
                2^nn! & \text{$n $ odd,}\\
                2^{n+1}n! & \text{$n$ even.}
            \end{cases}$\\[15pt]
            
            $2^nn!$ & $((C_n)_2,2Q(R^{\vee}),$\newline \mbox{}\hfill$2\omega_n^{\vee}+\al_n)$  &  $\ C_n$ & $1$ &  
            $\begin{cases}
                2^nn! & \text{$n$ odd,} \\ 
                2^{n+1}n! & \text{$n$ even.} 
            \end{cases}$ & 
            $\begin{cases}
                2^nn! & \text{$n $ odd,}\\
                2^{n+1}n! & \text{$n$ even.}
            \end{cases}$\\[10pt]  
    \noalign{\vskip 4pt}\hline\noalign{\vskip 4pt}
        \multicolumn{6}{l}{$D_n$ $n\ge3$ and $n\ne4,6$:}\\ 
            $1$ & $(D_n,Q(R),0)$& $\{1\}$ & $1$ & $1$ & $1$ \\ 
             
            $2$ & $(D_n^+,Q(R),0)$ & $\FC_2$ & $\V_1(D_n)-1$  & $1$ & $\V_1(D_n)-1$ \\             
                          
            $n!$ & $((D_n)_{n-1,n},Q(R),0)$  & $\FS_n$ & 
            $\begin{cases}
                2^{n-1} & \text{$n$ odd,}\\ 
                2^n & \text{$n$ even}.                
            \end{cases}$
            & $n!$ &
            $\begin{cases}
                2^{n-1}n! & \text{$n$ odd,}\\ 
                2^nn! & \text{$n$ even}.                
            \end{cases}$\\[15pt]

            $ 2^{n-1}n! $& $(\{1\},Q(R),0)$ & $D_n$ & $1$ &
            $\begin{cases}
                2^{n-1}n! & \text{$n$ odd,}\\
                2^nn! & \text{$n$ even.}
            \end{cases}$
            & $\begin{cases}
                2^{n-1}n! & \text{$n$ odd,}\\
                2^nn! & \text{$n$ even.}
            \end{cases}$\\[15pt]

            $ 2^{n-1}n! $& $((D_n)_{n-1,n},2Q(R)$
            \newline \mbox{}\hfill$+2\ZZ \omega_1,0)$ & $D_n$ & $1$ &
            $\begin{cases}
                2^{n-1}n! & \text{$n$ odd,}\\
                2^nn! & \text{$n$ even.}
            \end{cases}$
            & $\begin{cases}
                2^{n-1}n! & \text{$n$ odd,}\\
                2^nn! & \text{$n$ even.}
            \end{cases}$\\[15pt]
                         
            $ 2^{n-1}n! $& $((D_n)_{n-1,n},2Q(R)$
            \newline \mbox{}\hfill$+2\ZZ \omega_1,\al_{n-1})$ & $D_n$ & $1$ &
            $\begin{cases}
                2^{n-1}n! & \text{$n$ odd,}\\
                2^nn! & \text{$n$ even.}
            \end{cases}$
            & $\begin{cases}
                2^{n-1}n! & \text{$n$ odd,}\\
                2^nn! & \text{$n$ even.}
            \end{cases}$\\[10pt]           
            
        \multicolumn{6}{l}{$D_n$ ($n$ even):}\\ 
            $2^{n-2}n!$  & $(\{\pm1\},Q(R),0)$  & $D_n/\{\pm1\}$ & $0 $& N/A & $0$\\

            $2^{n-2}n!$  & $((D_n)_{n-1,n},2P(R),0)$  & $\FC_2\times D_n/\{\pm1\}$ & $0 $& N/A & $0$\\

            $2^{n-2}n!$  & $((D_n)_{n-1,n},2P(R),\al_{n-1})$  & $D_n/\{\pm1\}$ & $0$ & N/A & $0$\\
    \hline
    \hline
        \caption{The homomorphism tables of $(\widetilde{A}_n,A_n)$ when $n>3$ and $n\ne 5$; and $(\widetilde{B}_n,B_n)$, $(\widetilde{C}_n,C_n)$, and $(\widetilde{D}_n,D_n)$ when $n>2$ and $n\neq 4,6$. See Definition \ref{def: homomorphism table} for a description of an homomorphism table. Recall that $\V_k(G)$ is the number of ordered $k$-tuples of pairwise commuting elements in $G$ each of which has order at most 2.}
        \label{tab: counting Hom(W~,W) for A,B,C,D}
    \end{longtable}
    \endgroup
        
\section{The homomorphisms tables of $\Hom(\widetilde{W},W)$ for the exceptional affine reflection groups}

    The homomorphism tables and number of homomorphisms for the remaining irreducible affine reflection groups, $\widetilde{A}_2$, $\widetilde{A}_5$, $\widetilde{B}_2=\widetilde{C}_2$, $\widetilde{B}_4$, $\widetilde{B}_6$, $\widetilde{C}_4$, $\widetilde{C}_6$, $\widetilde{D}_4$, $\widetilde{D}_6$, $\widetilde{E}_6$, $\widetilde{E}_7$, $\widetilde{E}_8$, $\widetilde{F}_4$, and $\widetilde{G}_2$ are presented in Table \ref{tab: counting exceptions}. Note that $\widetilde{A}_3\cong \widetilde{D}_3$ so $\Hom(\widetilde{A}_3,A_3)$ is described in Table \ref{tab: counting Hom(W~,W) for A,B,C,D}. These tables were computed with the assistance of {\sc Magma} and {\rm\cite{random_endos}}.
    Descriptions of the normal subgroups and automorphism groups of Weyl groups can be found in {\rm\cite{maxwell_normal_subs}} and {\rm\cite[Ch. 2]{AutW_Franzsen}}, respectively. 
    
    \begingroup\fontsize{8}{10}\selectfont 
    \begin{longtable}{ m{0.4in} m{1.25in} m{0.55in} m{0.8in} m{0.67in} m{1.1in} | m{0.75in} }
            \label{tab: counting exceptions}

       $ |\widetilde{W}:H|$ & 
       $(K,L,(c(x_1),\dots,c(x_k))$\newline \mbox{}\hfill$\sim H $\hfill\mbox{} &  
       $\widetilde{W}/H$ &
       $|\{X\leqslant W:$ \newline\mbox{}\hfill$X\cong \widetilde{W}/H\}|$ & 
       $|\Aut(\widetilde{W}/H)|$ &
       $|\{\varphi\in \mathrm{Hom}(\widetilde{W},W):$\newline \mbox{}\hfill$\ker\varphi=H\}|$ & 
       $|\Hom(\widetilde{W},W)|$ \\ 
    \hline\hline
        \multicolumn{6}{l|}{$A_2$:}&\\ 
             

            
            $6$ & $(\{1\},Q(R),0)$  & $A_2$ & $1$ & $6$ & $6$ & $28$\ \\ 
    
            $6$ & $(\{1\},Q(R),0)$  & $A_2$ & $1$ & $6$ & $6$ &\\ [6pt]
    \hline 
        \multicolumn{6}{l|}{$A_5$:}&\\ 
             

            $720$ & $(\{1\},Q(R),0)$  & $A_5$ & $1$ & $1440$* & $1440*$ & $1516$ \\ [6pt]         
    \hline    
        \multicolumn{6}{l|}{$B_2=C_2$:}&\\ 
            $1$ & $(B_2,Q(R^{\vee}),0)$& $\{1\}$ & $1$ & $1$ & $1$ & $168$ \\ 
             
            $2$ & $(B_2^+,Q(R^{\vee}),0)$ & $\FC_2$ & $5$  & $1$ & $5$ &\\             
             
            $2$ & $((B_2)_1,Q(R^{\vee}),0) $ &  $\FC_2 $ & $5$ & $1$ & $5$ &\\
             
            $2$ & $((B_2)_2,Q(R^{\vee}),0) $ &  $\FC_2 $  & $5$ & $1$ & $5$ & \\
            $2$ & $(B_2,Q(R),(0,0))$& $\FC_2$ & $5$ & $1$ & $5$ &\\ 

            $2$ & $(B_2,Q(R),(0,\al_2))$& $\FC_2$ & $5$ & $1$ & $5$& \\ 

            $2$ & $(B_2,Q(R),(\al_2,0))$& $\FC_2$ & $5$ & $1$ & $5$ &\\ 

            $2$ & $(B_2,Q(R),(\al_2,\al_2))$& $\FC_2$ & $5$ & $1$ & $5$ &\\ 
             
            $4$ & $(\{\pm1\} ,Q(R^{\vee}),0)$ & $\FC_2^2 $   & $2$& $6$ & $12$& \\
            $4$ & $((B_2)_1 ,Q(R),0)$ & $\FC_2^2 $   & $2$& $6$ & $12$ &\\

            $4$ & $((B_2)_1 ,Q(R),\al_2)$ & $\FC_2^2 $   & $2$& $6$ & $12$ &  \\

            $4$ & $((B_2)_2,Q(R),0)$ & $\FC_2^2 $   & $2$& $6$ & $12$ &\\
            $4$ & $((B_2)_2,Q(R),\al_2)$ & $\FC_2^2 $   & $2$& $6$ & $12$ &\\
            $4$ & $(B_2^+ ,Q(R),0)$ & $\FC_2^2 $   & $2$& $6$ & $12$ &\\
            
            $4$ & $(B_2^+ ,Q(R),\al_2)$ & $\FC_2^2 $   & $2$& $6$ & $12$& \\

            $8$ & $(\{1\} ,Q(R^{\vee}),0)$ & $B_2$   & $1$ & $8$ & $8$ &\\

            $8$ & $(\{\pm1\} ,Q(R),0)$ & $B_2$   & $1$ & $8$ & $8$ &\\
            $8$ & $(\{\pm1\} ,Q(R),\alpha_2)$ & $\FC_2^3$   & $0$ & $8$ & $0$ &\\
            $8$ & $((B_2)_2 ,2Q(R^{\vee}),0)$ & $B_2$   & $1$ & $8$ & $8$ &\\

            $8$ & $((B_2)_2  ,2Q(R^{\vee}),\al_2)$ & $B_2$   & $1$ & $8$ & $8$ &\\

            $8$ & $((B_2)_2  ,2Q(R^{\vee}),2\omega_2^{\vee})$ & $B_2$   & $1$ & $8$ & $8$ &\\

            $8$ & $((B_2)_2 ,2Q(R^{\vee}),$\newline \mbox{}\hfill$2\omega_2^{\vee}+\al_2)$ & $B_2$   & $1$ & $8$ & $8$ &\\[6pt]        
         \hline
        \multicolumn{6}{l|}{$B_4$:}&\\ 
            $6$ & $(N,Q(R^{\vee}),(0,0))$ & $\FS_3$ & $64$& $6$ & $384$ & $14176$\\

            $12$ & $(N^+,Q(R^{\vee}),0)$ & $\FC_2\times\FS_3$ & $96$& $12$ & $1152$ & \\

            $24$ & $((B_4)_2,Q(R^{\vee}),0)$ & $\FS_4$ & $2^5$*& $4!$ & $2^54!$* &  \\

            $24$ & $(N,2P(R),(0,0))$ & $\FS_4$ & $2^5$& $4!$ & $2^54!$ &  \\

            $24$ & $(N,2P(R),(0,\al_3))$ & $\FS_4$ & $2^5$& $4!$ & $2^54!$ &  \\

            $48$ & $((B_4)_2^+,Q(R^{\vee}),0)$ & $\FC_2\times\FS_4$ & $2^5$*& $2\cdot5!$ & $2^64!$* &  \\

            $48$ & $(N^+,2P(R),(0,0))$ & $\FC_2\times \FS_4$ & $2^5$& $2\cdot5!$ & $2^64!$ &  \\

            $48$ & $(N^+,2P(R),(0,\al_3))$ & $\FC_2\times\FS_4$ & $2^5$ & $2\cdot5!$ & $2^64!$ &  \\[6pt]
    \hline
        \multicolumn{6}{l|}{$B_6$:}&\\ 
            $6!$ & $((B_6)_2,Q(R^{\vee}),0)$  & $\mathfrak{S}_6$ & $2^6 $& $2\cdot 6!$* & $2^76!$*& $660736$\\

            $2 \cdot 6!$& $((B_6)_2^+,Q(R^{\vee}),0)$ & $\FC_2\times \mathfrak{S}_6$  & $2^{5}$& $2^2\cdot 6! $*& $2^76!$* & \\[6pt]  
    \hline
        \multicolumn{6}{l|}{$C_4$:}&\\ 
            $6$ & $(N,Q(R^{\vee}),(0,0))$ & $\FS_3$ & $64$& $6$ & $384$ & $54112$\\

            $12$ & $(N,Q(R)),(0,0))$ & $\FC_2\times\FS_3$ & $96$& $12$ & $1152$ & \\
                                    
            $12$ & $(N,Q(R),(\al_4,0))$ & $\FC_2\times\FS_3$ & $96$& $12$ & $1152$ & \\

            $12$ & $(N^+,Q(R^{\vee}),0)$ & $\FC_2\times\FS_3$ & $96$& $12$ & $1152$ & \\

            $24$ & $(N^+,Q(R),0)$ & $\FC_2^2\times\FS_3$ & $16$& $144$ & $2304$ & \\

            $24$ & $((C_4)_2,Q(R^{\vee}),0)$ & $\FS_4$ & $32$*& $24$ & $768$* & \\

            $48$ & $((C_4)_2^+,Q(R^{\vee}),0)$ & $\FC_2\times\FS_4$ & $32$*& $48$ & $1536$* & \\

            $48$ & $((C_4)_2,Q(R),0)$ & $\FC_2\times \FS_4$ & $32$*& $48$ & $1536$* & \\

            $48$ & $((C_4)_2,Q(R),\al_4)$ & $\FC_2\times\FS_4$ & $32$*& $48$ & $1536$* &  \\

            $96$ & $((C_4)_2^+,Q(R),0)$ & $\FC_2^2\times\FS_4$ & $4$*& $576$ & $2304$* & \\[6pt]
    \hline
        \multicolumn{6}{l|}{$C_6$:}&\\ 
            $6!$ & $((C_6)_2,Q(R^{\vee}),0)$  & $\mathfrak{S}_6$ & $2^6 $& $2\cdot 6!$* & $2^76!$*& $17499904$\\

            $2 \cdot 6!$& $((C_6)_2^+,Q(R^{\vee}),0)$ & $\FC_2\times \mathfrak{S}_6$  & $2^{5}$& $2^2\cdot 6! $*& $2^76!$* & \\

            $2 \cdot 6!$& $((C_6)_2,Q(R),0)$ & $\FC_2\times \mathfrak{S}_6$  & $2^{5}$& $2^2\cdot 6! $*& $2^76!$*& \\

            $2 \cdot 6!$& $((C_6)_2,Q(R),\al_6)$ & $\FC_2\times \mathfrak{S}_6$  & $2^{5}$& $2^2\cdot 6! $*& $2^76!*$ & \\[6pt]    
    \hline
        \multicolumn{6}{l|}{$D_4$:}&\\ 
            $6$ & $(N,Q(R),(0,0))$ & $\FS_3$ & $32$& $6$ & $192$ & $12332$\\

            $24$ & $((D_4)_{3,4},Q(R),0)$  & $\mathfrak{S}_4$ & $24$* & $24$ & $576$* &  \\ 

            $24$ & $((D_4)_{1,4},Q(R),0)$  & $\mathfrak{S}_4$ & $24$ & $24$ & $576$ & \\ 

            $24$ & $((D_4)_{1,3},Q(R),0)$  & $\mathfrak{S}_4$ & $24$ & $24$ & $576$ & \\ 

            $24$ & $(N,2P(R),(0,0))$  & $\mathfrak{S}_4$ & $24$ & $24$ & $576$ & \\ 

            $24$ & $(N,2P(R),(\al_1,0))$  & $\mathfrak{S}_4$ & $24$ & $24$ & $576$ & \\ 

            $24$ & $(N,2P(R),(0,\al_1))$  & $\mathfrak{S}_4$ & $24$ & $24$ & $576$ & \\ 

            $24$ & $(N,2P(R),(\al_1,\al_1))$  & $\mathfrak{S}_4$ & $24$ & $24$ & $576$ & \\ 

            $96$ & $((D_4)_{1,4},2P(R),0)$  & $D_4/\{\pm1\}$ & $0$ & $576$ & $0$ & \\ 

            $96$ & $((D_4)_{1,4},2P(R),\al_)$  & $D_4/\{\pm1\}$ & $0$ & $576$ & $0$ & \\ 

            $96$ & $((D_4)_{1,3},2P(R),0)$  & $D_4/\{\pm1\}$ & $0$ & $576$ & $0$ & \\ 

            $96$ & $((D_4)_{1,3},2P(R),\al_1)$  & $D_4/\{\pm1\}$ & $0$ & $576$ & $0$ & \\ 

            $1$ & $(\{1\},Q(R),0)$ & $D_4$ & $1$ & $1152$* & $1152$* &  \\
            
            $192$& $((D_4)_{3,4},2Q(R)$\newline \mbox{}\hfill$+2\ZZ \omega_1,0)$ & $D_4$ & $1$ & $1152$* & $1152$* & \\

            $192$& $((D_4)_{3,4},2Q(R)$\newline \mbox{}\hfill$+2\ZZ \omega_1,\al_{3})$ & $D_4$ & $1$ & $1152$* & $1152$* &  \\

            $192$& $((D_4)_{1,4},2Q(R)$\newline \mbox{}\hfill$+2\ZZ \omega_1,0)$ & $D_4$ & $1$ & $1152$ & $1152$ &\\
            
            $192$& $((D_4)_{1,4},2Q(R)$\newline \mbox{}\hfill$+2\ZZ \omega_1,\al_{1})$ & $D_4$ & $1$ & $1152$ & $1152$ &\\           

            $192$& $((D_4)_{1,3},2Q(R)$\newline \mbox{}\hfill$+2\ZZ \omega_1,0)$ & $D_4$ & $1$ & $1152$ & $1152$ &\\

            $192$& $((D_4)_{1,3},2Q(R)$\newline \mbox{}\hfill$+2\ZZ \omega_1,\al_{1})$ & $D_4$ & $1$ & $1152$ & $1152$ &                        \\[6pt]
    \hline
        \multicolumn{6}{l|}{$D_6$:}&\\ 
            $6!$ & $((D_6)_{5,6},Q(R),0)$  & $\mathfrak{S}_6$ & $2^6 $& $2\cdot 6!$* & $2^76!$* & $231152$\\[6pt]         
    \hline
        \multicolumn{6}{l|}{$E_6$:}&\\ 
            $1$ & $(E_6,Q(R),0)$& $\{1\}$ & $1$ & $1$ & $1$ & $52732$ \\ 
             
            $2$ & $(E_6^+,Q(R),0)$ & $\FC_2$ & $891$  & $1$ & $891$ &  \\   

            $51840$ & $(\{1\},Q(R),0)$ & $E_6$ & $1$ & $51840$ & $51840$ \\ [6pt]
    \hline
        \multicolumn{6}{l|}{$E_7$:}&\\ 
            $1$ & $(E_7,Q(R),0)$& $\{1\}$ & $1$ & $1$ & $1$ & $2913248$ \\ 
             
            $2$ & $(E_7^+,Q(R),0)$ & $\FC_2$ & $10207$  & $1$ & $10207$ &  \\  

            $1451520$ & $(Z(E_7),Q(R),0)$ & $\Inn(E_7)$& $1$ & $1451520$  & $1451520$ &  \\
                    
            $2903040$ & $(\{1\},Q(R),0)$  & $E_7$ & $1$ & $1451520$ & $1451520$ & \\ [6pt]
    \hline
        \multicolumn{6}{l|}{$E_8$:}&\\ 
            $1$ & $(E_8,Q(R),0)$& $\{1\}$ & $1$ & $1$ & $1$ & $696929552$ \\ 
            $2$ & $(E_8^+,Q(R),0)$ & $\FC_2$ & $199951$  & $1$ & $199951$ &  \\ 
            $348364800$ & $(\{\pm1\},Q(R),0)$ & $\{\pm1\}$ & $0$ & $\mathrm{N/A}$ & $0$ & \\
            $696729600$ & $(\{1\},Q(R),0)$  & $E_8$ & $1$ & $696729600$ & $696729600$ & \\ [6pt]
    \hline
        \multicolumn{6}{l|}{$F_4$:}& \\ 
            $1$ & $(F_4,Q(R^{\vee}),0)$& $\{1\}$ & $1$& $1$& $1$ & $40096$ \\ 
            $2$ & $(F_4^+,Q(R^{\vee}),0)$ & $\FC_2$  & $139$ & $1$ & $139$ & \\
            $2$ & $(\prescript{}{+}{F_4},Q(R^{\vee}),0)$ & $\FC_2$  & $139$ & $1$ & $139$ & \\
            $2$ & $((F_4)_+,Q(R^{\vee}),0)$ & $\FC_2$ & $139$ & $1$ & $139$ & \\
            $4$ & $(\prescript{}{+}({F_4})_+,Q(R^{\vee}),0)$ & $\FC_2^2$ & $597$ & $6$ & $3582$ & \\ 
            $6$ & $((F_4)_1,Q(R^{\vee}),0)$ & $\mathfrak{S}_3$ & $352$ & $6$ & $2112$ &  \\  
            $6$ & $((F_4)_2,Q(R^{\vee}),0)$ & $\mathfrak{S}_3$ & $352$ & $6$ & $2112$ & \\  
            $12$ & $((F_4)_1^+,Q(R^{\vee}),0)$ & $\FC_2\times \mathfrak{S}_3$ & $560$ & $12$ & $6720$ & \\  
            $12$ & $((F_4)_2^+,Q(R^{\vee}),0)$ & $\FC_2\times \mathfrak{S}_3$ & $560$ & $12$ & $6720$ & \\ 
            $24$ & $((F_4)_2,2Q(R),0)$ & $\mathfrak{S}_4$ & $96$ & $24$ & $2304$ & \\  
            $36$ & $(N,Q(R^{\vee}),0)$ & $(\mathfrak{S}_3)^2$ & $64$ & $72$ & $4608$ & \\ 
            $48$ & $((F_4)_2^+,2Q(R),0)$ & $\FC_2\times \mathfrak{S}_4$ & $144$ & $48$ & $6912$ & \\ 
            $144$ & $(N,2Q(R),0)$ & $\FS_4\times\mathfrak{S}_3$ & $0$ & $144$ & $0$ & \\ 
            $576$ & $(\{\pm 1\},Q(R^{\vee}),0)$ & $\Inn(F_4)$ & $0$ & $1152$ & $0$ & \\ 
            $1152$ & $(\{1\},Q(R^{\vee}),0)$ & $F_4$ & $1$ & $4608$ & $4608$ & \\ [6pt]
    \hline
        \multicolumn{6}{l|}{$G_2$:}&\\ 
            $1$ & $(G_2,Q(R^{\vee}),0)$& $\{1\}$ & $1$& $1$& $1$ & $88$ \\ 
            $2$ & $(G_2^+,Q(R^{\vee}),0)$ & $\FC_2$  & $7$ & $1$ & $7$ & \\
            $2$ & $((G_2)_1,Q(R^{\vee}),0)$ & $\FC_2$  & $7$ & $1$ & $7$ & \\
            $2$ & $((G_2)_2,Q(R^{\vee}),0)$ & $\FC_2$ & $7$ & $1$ & $7$ & \\
            $4$ & $(\lb (s_1s_2)^2\rb,Q(R^{\vee}),0)$ & $\FC_2^2$ & $3$ & $6$ & $18$ &  \\ 
            $6$ & $(\{\pm1\},Q(R^{\vee}),0)$ & $\mathfrak{S}_3$ & $2$ & $6$ & $12$ & \\  
            $6$ & $((G_2)_1,Q(R),0)$ & $\mathfrak{S}_3$ & $2$ & $6$ & $12$ & \\  
            $12$ & $(\{1\},Q(R^{\vee}),0)$ & $G_2$ & $1$ & $12$ & $12$ & \\  
            $12$ & $(\lb (s_1s_2)^2\rb,Q(R),0)$ & $G_2$ & $1$ & $12$ & $12$ & \\ 
    \caption{The homomorphism tables of $(\widetilde{W},W)$ where $\widetilde{W}$ is an exceptional irreducible affine reflection group. See Definition \ref{def: homomorphism table} for a description of an homomorphism table. For types $A_2$, $A_5$, $B_4$, $B_6$, $C_4$, $C_6$, $D_4$ and $D_6$, an asterisk * in a row indicates that this row \textit{replaces} the corresponding row in Table \ref{tab: counting Hom(W~,W) for A,B,C,D}, and asterisked entries are those that differ from those in Table \ref{tab: counting Hom(W~,W) for A,B,C,D}. Rows without asterisks should be taken in \textit{addition} to those in Table \ref{tab: counting Hom(W~,W) for A,B,C,D}. For type $C_2= B_2$, the complete homomorphism table is given.}
    \end{longtable}
    \endgroup

\section{Random homomorphisms}
\label{sec: random homs}

    In this section we describe properties of random homomorphisms $\widetilde{C}_n\to C_n$. Similar proofs describe random homomorphisms $\widetilde{W}_n\to W_n$ where $W_n=A_n$, $B_n$, or $D_n$. We are interested in examining the behaviour of $\Hom(\widetilde{C}_n,C_n)$ as $n\to \infty$ so, to avoid the exceptional homomorphisms that occur when $n=2,4$ or $n=6$, we assume for the remainder of this section that $n\geq 7.$
    Let $X_n$ be the order of the image of an homomorphism drawn uniformly at random from $\Hom(\widetilde{C}_n,C_n)$. Then the probability that a random homomorphism $\varphi\in \Hom(\widetilde{C}_n,C_n)$ satisfies $|\im(\varphi)|=k$ is 
    \begin{equation}
        \PP\big[X_n=k\big]=
        \frac{1}{|\Hom(\widetilde{C}_n,C_n)|}\sum_{\substack{K\unlhd \widetilde{C}_n \\ |\widetilde{C}_n:K|=k }} \big|\{\varphi\in \Hom(\widetilde{C}_n,C_n):\ker\varphi=K\}\big|     .
    \label{equ: random variable general}
    \end{equation} 
    The homomorphism table of $\Hom(\widetilde{C}_n,C_n)$ (shown in Table \ref{tab: counting Hom(W~,W) for A,B,C,D}) lists, in order of increasing index, the normal subgroups $K$ of $\widetilde{C}_n$ along with $\big|\{\varphi\in \Hom(\widetilde{C}_n,C_n):\ker\varphi=K\}\big|$, the number of homomorphisms  with kernel $K$, and so this table allows us to compute $\PP[X_n=k]$. 
    \begin{equation}
     \PP\big[X_n=k\big]\cdot \big| \Hom(\widetilde{C}_n,C_n)\big|=
    \begin{cases}
        1 & k=1,\\
        7(\mathcal{V}_1(C_n)-1) & k=2,\\
        7(\mathcal{V}_2(C_n)-3\mathcal{V}_1(C_n)+2) & k=4,\\
        \mathcal{V}_3(C_n)-7\mathcal{V}_2(C_n)+14\mathcal{V}_1(C_n)-8 & k=8,\\
        2^nn! & k=n! ,\\
        3\cdot 2^nn! & k=2\cdot n!, \\ 
        3\cdot 2^nn! & \text{$n$ is odd and }k=2^{n-1}n!,\\
        9\cdot 2^nn! & \text{$n$ is odd and $k=2^nn!$},\\
        3\cdot 2^{n+2}n! & \text{$n$ is even and $k=2^nn!$},\\ 
        0 & \text{else.}
    \end{cases}
    \label{equ: random variable Cn}
    \end{equation}

    We first show that the sequences $(\mathcal{V}_d(C_n))_{n\in\NN} $ where $d\in \{1,2,3\}$ grow more slowly that $n!$. Then, since the number of subgroups of $C_n$ isomorphic to $\FC_2^d$ is a linear combination of $\V_1(C_n),\dots, \V_d(C_n)$ and $|\Hom(\widetilde{C}_n,C_n)|>n!$, this will allow us to conclude that homomorphisms $\varphi:\widetilde{C}_n\to C_n$ asymptotically almost surely have image-order greater than 8.
    
    \begin{lemma} For $d\in \{1,2,3\}$, 
    \label{lem: V_k(C_n)/n! tends to 0}
        $$ \lim_{n\to\infty} \frac{\mathcal{V}_d(C_n)}{n!}\to 0. $$
    \end{lemma}
    \begin{proof}
        Since $\mathcal{V}_d(G)$ counts the number of tuples $(x_1,\dots,x_d)$ of pairwise commuting elements $x_i\in G$ with $x_i^2=1$, clearly $\mathcal{V}_d(G)\leq \mathcal{V}_{d+1}(G)$. Hence it suffices to prove that 
        $$\frac{\mathcal{V}_3(C_n)}{n!}= \sum_{k=0}^{\lfloor \frac{n}{8}\rfloor} \;
            \sum_{l=0}^{\lfloor \frac{n}{4}\rfloor-2k}\; 
            \sum_{m=0}^{\lfloor \frac{n}{2}\rfloor-4k-2l} 
            \frac{2^{3n}7^{l+m}}{2^{20k+10l+4m} \ k!\ l!\ m!\ (n-8k-4l-2m)!}  \xrightarrow[]{n\to \infty} 0. $$
        
        At least one of $k$, $l$, $m$, or $n-8k-4l-2m$ is at least $\frac{n}{15}$, so the product of their factorials is at least $\lfloor \frac{n}{15}\rfloor !$ and hence
        $$ \frac{\mathcal{V}_3(C_n)}{n!}\leq  
        \sum_{k=0}^{\lfloor \frac{n}{8}\rfloor} \;
            \sum_{l=0}^{\lfloor \frac{n}{4}\rfloor-2k} \;
            \sum_{m=0}^{\lfloor \frac{n}{2}\rfloor-4k-2l}\frac{2^{3n}7^{l+m}}{\lfloor \frac{n}{15}\rfloor !}\leq \frac{2^{3n}7^{n} n^3}{\lfloor \frac{n}{15}\rfloor !}\xrightarrow[]{n\to \infty}0.
        $$
    \end{proof}

    Lemma \ref{lem: V_k(C_n)/n! tends to 0} allows us to study the asymptotic behaviour of homomorphisms. 
    
    \begin{theorem}
    \label{thm: asymptotic image size}
    As $n$ tends to infinity, a random homomorphism  $\widetilde{C}_n\to C_n$ almost surely has large image and small kernel.
    \begin{enumerate}[label=(\alph*)]
        \item The image-order of a homomorphism $\widetilde{C}_n\to C_n$ is asymptotically almost surely between $n!$ and $2^nn!=|C_n|$.
        $$ \lim_{n\to \infty} \PP\big
        [ n! \leq X_n\leq 2^nn!\big]
         = 1. $$
         \item The number of homomorphisms $\widetilde{C}_n\to C_n$ is asymptotically $2^{n+4}n!$,
         $$ \big|\Hom(\widetilde{C}_n,C_n)\big|\sim 2^{n+4}n!. $$
    \end{enumerate}
    \end{theorem}
    \begin{proof}
        By Theorem \ref{thm: main theorem}, we know that $\big|\Hom(\widetilde{C}_n,C_n)\big|=2^{n+4}n!+\mathcal{V}_3(C_n)$. Moreover, from Table \ref{tab: counting Hom(W~,W) for A,B,C,D} and Lemma \ref{lem: V_k(C_n)/n! tends to 0}, we find that
        $$ \PP\big[ 1 \leq X_n <n! \big]=\PP\big[ 1 \leq X_n \leq8 \big] = \frac{\mathcal{V}_3(C_n)}{|\Hom(\widetilde{C}_n,C_n)|}\leq \frac{\V_3(C_n)}{n!} \xrightarrow[]{n\to \infty} 0.$$
        Since the order of the image of a homomorphism $\widetilde{C}_n\to C_n$ lies between $1$ and $2^nn!=|C_n|$, this proves part (a).
    
        For part (b), observe that
        $$ \lim_{n\to\infty}\frac{\big|\Hom(\widetilde{C}_n,C_n)\big|}{2^{n+4}n!}=\lim_{n\to\infty}\frac{2^{n+4}n!+\mathcal{V}_3(C_n)}{2^{n+4}n!}=1  .$$
    \end{proof}

    Equation \eqref{equ: random variable Cn} lists $\PP[X_n=k]$ for any $k\in \NN$, so we can compute the expected value of the image-order of a homomorphism $\widetilde{C}_n\to C_n$. 

    \begin{proof}[Proof of Theorem \ref{thm: expected image order A,B,C,D}]
        The order of the image of a homomorphism $\widetilde{C}_n\to C_n$ lies between $1$ and $2^nn!$. 
        For odd $n$, Equation \eqref{equ: random variable Cn} implies that the expectation is
        $$\EE[X_{n}]= \frac{1}{|\Hom(\widetilde{C}_n,C_n)|}\Big( -21+41\mathcal
        V_1(C_n)-28\mathcal{V}_2(C_n)+8\mathcal{V}_3(C_n)+2^n(n!)^2(7+21\cdot 2^{n-1}) \Big).
        $$
        Using Lemma \ref{lem: V_k(C_n)/n! tends to 0} and Theorem \ref{thm: asymptotic image size}, we find that 
        $$ 
        \EE[X_{2n+1}]\sim \frac{21\cdot 2^{2n}(2n+1)!}{2^{2n+5}(2n+1)!}=\frac{21}{32}|C_{2n+1}|.
        $$
    \end{proof}

    Since homomorphisms $\widetilde{C}_n\to C_n$ with image-order $2^nn!$ are surjections, by computing $\PP[X_n=2^nn!]$ we find the probability that a random $\varphi\in \Hom(\widetilde{C}_n,C_n)$ is surjective.

    \begin{proof}[Proof of Theorem \ref{thm: prob surjective hom ABCD}]
        From Table \ref{tab: counting Hom(W~,W) for A,B,C,D}, we see that there are $3\cdot 2^{2n+2}(2n)!$ surjections $\widetilde{C}_{2n}\twoheadrightarrow C_{2n}$ and by Theorem \ref{thm: asymptotic image size} part (b), the number of homomorphisms $\widetilde{C}_{2n}\to C_{2n}$ is asymptotically $2^{2n+4}(2n)!$. Hence,
        $$ \PP\big[ X_{2n}=2^{2n}(2n)!]= \frac{3\cdot 2^{2n+2}(2n)!}{|\Hom(\widetilde{C}_{2n},C_{2n})|}\sim \frac{3\cdot 2^{2n+2}(2n)!}{2^{2n+4}(2n)!}=\frac{3}{4}. $$
    \end{proof}

    We can also examine how generic $\pi_{\text{standard}}$ is in the set of surjective homomorphisms $\widetilde{C}_n\twoheadrightarrow C_n$.

    \begin{proof}[Proof of Theorem \ref{thm: prob surjective factors thru standard ABCD}]
        From Table \ref{tab: counting Hom(W~,W) for A,B,C,D} we see that there are $3\cdot 2^{2n+2}(2n)!$ surjective homomorphisms $\widetilde{C}_{2n}\twoheadrightarrow C_{2n}$, $2^{2n+1}(2n)!$ factor through $\pi_{\text{standard}}$.
    \end{proof}

    Similarly, we investigate how generic homomorphisms which factor through $\pi_{\text{standard}}$ are in $\Hom(\widetilde{C}_n,C_n).$

    \begin{proof}[Proof of Theorem \ref{thm: prob factors thru standard ABCD}]
        A homomorphism $\varphi:\widetilde{C}_n\to C_n$ factors through $\pi_{\mathrm{standard}}$ if and only if its kernel $\ker\varphi\sim (K,L,(c(x_1),\dots,c(x_k)))$ has $L=Q(R^{\vee})$. 
        $$ \PP_{\varphi\in \Hom(\widetilde{C}_n,C_n)}\big[\varphi \text{ factors through }\pi_{\text{standard}} \big]=\frac{\mathcal{V}_2(C_n)+2^{n+2}n!}{|\Hom(\widetilde{C}_n,C_n)|}\sim \frac{2^{n+2}n!}{2^{n+4}n!}=\frac{1}{4}. $$
    \end{proof}

\section{Describing homomorphisms $\widetilde{W}\to W$ explicitly}
\label{sec: explicit homomorphisms}

    In this section we explicitly write down the homomorphisms $\widetilde{W}_n\to W_n$ on the standard generating set $s_0,\dots,s_n$ of $\widetilde{W}$ and, when $W_n$ is $B_n$, $C_n$, or $D_n$, on an arbitrary element of $W_n$. 
    Our aim is to provide a unified and simple description of these homomorphisms using our understanding of the normal subgroups $N\unlhd \widetilde{W}_n$ and the equivalence classes of subgroups of $W_n$ isomorphic to $\widetilde{W}_n/N$ under the action of $\Aut(W_n)$. We attempt to list homomorphisms with minimal redundancy, but prioritise simplicity in the statement where these conflict. 
    To avoid the exceptional cases, in this section we assume that $n\geq 4$ when $W_n=A_n$, and $n\geq 3$ and $n\neq 4$ when $W_n=B_n$, $C_n$, or $D_n$.

\subsection{Types $A$, $E$, $F$ and $G$}

    When $W_n$ is $A_n$, $E_6$, $E_7$, or $E_8$, we see from Tables \ref{tab: counting Hom(W~,W) for A,B,C,D} and \ref{tab: counting exceptions} that every homomorphism $\widetilde{W}_n\to W_n$ factors through $\pi_{\mathrm{standard}}$. Hence to describe $\Hom(\widetilde{W}_n,W_n)$, it suffices to describe $\Hom(W_n,W_n)$. If $W_n=A_n$, $E_6$ or $E_7$, then every endomorphism of $W_n$ either has image-order at most 2, or is an automorphism. For $E_7$, there is also an endomorphism whose kernel is $Z(E_7)$. 
    
    \begin{theorem}
    \label{thm: explicit homs types A and E}
        Suppose $W_n=A_n$ with $n\geq 4$, or $W_n=E_6$, $E_7$, $E_8$. Every homomorphism $\varphi:\widetilde{W}_n\to W_n$ is the composition of $\pi_{\mathrm{standard}}$ with either: $ \mathrm{sgn}_{{W}_n}^w$ for some $w\in W_n$ with $w^2=1$; an automorphism of $W_n$; or, if $W_n=E_7$, the endomorphism of $E_7$ defined by $s_i\mapsto z_0s_i$ for $i\in \{1,\dots,n\}$ composed with an automorphism of $E_7$. 
    \end{theorem}

    Suppose $W_n$ is $F_4$ or $G_2$. From Table \ref{tab: counting exceptions}, we know what subgroups of $N\unlhd \widetilde{W}_n$ are kernels of homomorphisms $\widetilde{W}_n\to W_n$, and what the isomorphism types of $\widetilde{W}_n/N$ are. 
    Using {\sc Magma}, we find the subgroups of $W_n$ isomorphic to $\widetilde{W}_n/N$ and classify them up to equivalence under the action of $\Aut(W_n)$. For each $N\unlhd \widetilde{W}_n$ and $H\leqslant W_n$ with $H\cong \widetilde{W}_n/N$, we construct a homomorphism explicitly on the standard generators of $\widetilde{W}$ which has kernel $N$ and image $H$. Finally, we endeavour to compose our explicitly-constructed homomorphisms with automorphisms of $W_n$ so that the full list can be stated relatively simply.   

    Using Table \ref{tab: counting exceptions},  we see that there are $13$ normal subgroups of $\widetilde{F}_4$ which appear as kernels of homomorphisms $\widetilde{F}_4\to F_4$, and 9 different isomorphism types of images of homomorphisms $\widetilde{F}_4\to F_4$. 
    For example, there are two normal subgroups of $N\unlhd \widetilde{F}_4$ for which $\widetilde{F}_4/N$ is isomorphic to $\FS_3$, and two for which $\widetilde{F}_4/N$ is isomorphic to $\FC_2\times \FS_3$. Up to $\Aut(F_4)$, there are three subgroups of $F_4$ isomorphic to $\FS_3$, and four isomorphic to $\FC_2\times \FS_3$. These 14 (up to postcomposition with $\Aut(F_4)$) homomorphisms $\widetilde{F}_4\to F_4$ appear, with some redundancy, in the 16 (up to postcomposition with $\Aut(F_4)$) homomorphisms in part $(b)$ of Theorem \ref{thm: explicit homs type F}. The sum $\eps_3+\eps_4$ determines whether the image is $\FS_3$ or $\FC_2\times \FS_3$, the  permutation $\sigma \in \lb (13)(24)\rb$ dictates the specific kernel of the homomorphism, and the coefficients $\eps_1,\eps_2, \eps_3$, and $\eps_4$ decide the automorphism-class of the image.
    
    \begin{theorem}
    \label{thm: explicit homs type F}
        Suppose $\varphi:\widetilde{F}_4\to F_4$ is a homomorphism. Then one of the following is true.
        \begin{enumerate}[label=(\alph*)]
            \item The image of $\varphi$ satisfies $|\im(\varphi)|\leq4$ and there exist commuting $w_1,w_2\in F_4$ with $w_i^2=1$ so that 
            $$ \varphi(s_0)=\varphi(s_1)=\varphi(s_2)=w_1 \quad\text{and}\quad \varphi(s_3)=\varphi(s_4)=w_2. $$
            \item There exist $\epsilon_1,\eps_2,\eps_3,\eps_4\in \{0,1\}$ satisfying $\eps_1+\eps_2+\eps_3\leq 1$ and $\sigma\in \lb (13)(24)\rb$ so that, up to postcomposition with an automorphism of $F_4$, $\varphi$ is
            $$ \varphi(s_0)=r_{\widetilde{\alpha}}\quad 
            \varphi(s_{\sigma(1)})=s_3r_{\widetilde{\alpha}}^{\eps_1} \quad  
            \varphi(s_{\sigma(2)})=s_4r_{\widetilde{\alpha}}^{\eps_1}s_1^{\eps_2} \quad 
            \varphi(s_{\sigma(3)})=\varphi(s_{\sigma(4)})=s_1^{\eps_3}z_0^{\eps_4} .
            $$
            \item Up to postcomposition with an automorphism of $F_4$, $\varphi$ is
            $$ \varphi(s_0)=\varphi(s_2)=r_{\widetilde{\alpha}}\quad
            \varphi(s_i)=s_i \quad \text{for }i\in \{1,3,4\}.
            $$
            \item There exist $\eps_1,\eps_2,\eps_3\in \{0,1\}$ so that, up to postcomposition with an automorphism of $F_4$, $\varphi$ is
            $$ \varphi(s_0)=(s_0s_4)^{\eps_1}r_{\widetilde{\alpha}} \quad 
            \varphi(s_1)=(s_0s_4)^{\eps_1} s_1 \quad 
            \varphi(s_2)=(s_0s_4)^{\eps_1} s_2 \quad 
            \varphi(s_3)=\varphi(s_4)=z_0^{\eps_2}s_4^{\eps_3}$$
        \end{enumerate}
        Moreover,in cases $(a)$-$(c)$, $\varphi$ factors through $\pi_{\mathrm{standard}}$.
    \end{theorem}

    The homomorphism table of $G_2$, which appears in Table \ref{tab: counting exceptions}, shows us that $\widetilde{G}_2$ has 9 normal subgroups which appear as kernels of homomorphisms $\widetilde{G}_2\to G_2$, and 5 isomorphism classes of images of these homomorphisms. Examination of these subgroups yields Theorem \ref{thm: explicit homs type G}.

    \begin{theorem}
        \label{thm: explicit homs type G}
         Suppose $\varphi:\widetilde{G}_2\to G_2$ is a homomorphism. Then one of the following is true.
        \begin{enumerate}[label=(\alph*)]
            \item The image of $\varphi$ satisfies $|\im(\varphi)|\leq4$ and there exist commuting $w_1,w_2\in G_2$ with $w_i^2=1$ so that 
            $$ \varphi(s_1)=w_1 \quad\text{and}\quad \varphi(s_0)=\varphi(s_2)=w_2; $$
            \item The image of $\varphi$ satisfies $6\leq|\im(\varphi)|\leq12$ and, up to postcomposition with an automorphism of $G_2$, $\varphi$ is 
            $$ \varphi(s_0)=r_{\widetilde{\alpha}} \qquad \varphi(s_1)\in \{1,r_{\widetilde{\alpha}},z_0,s_1\} \qquad \varphi(s_2)=s_2.$$
        \end{enumerate}
        Moreover, $\varphi$ factors through $\pi_{\mathrm{standard}}$, except in case $(b)$ when $\varphi(s_1)\in \{1,z_0\}$.
    \end{theorem}

\subsection{Types $B$, $C$, and $D$}

    Suppose $W_n$ is $B_n,$ $C_n$, or $D_n$ and let $N$ be the coordinate inversion subgroup of $W_n$. Then $N$ is the unique index $n!$ normal subgroup of $W_n$, and is itself generated by reflections. Specifically,
    $$ 
    N=\begin{cases}
        (W_n)_2=\lb r_1,\dots,r_n\rb\cong \FC_2^n & W_n=B_n,C_n,\\ 
        (D_n)_{n-1,n}=\lb r_{-1,1,2},\dots,r_{-1,n-1,n}\rb\cong \FC_2^{n-1} &W_n=D_n.
    \end{cases}
    $$    
    Any element of $\widetilde{W}_n$ can be uniquely written in the form $(v,x,\sigma)$ where $v\in Q(R^{\vee})\cong \ZZ^n$, $x\in N$, and $\sigma\in S_+\cong \FS_n$ since $\widetilde{W}_n$ decomposes as
    $$ (v,x,\sigma)\in Q(R^{\vee})\rtimes (N\rtimes S_+)\cong \widetilde{W}_n .$$
    Using this decomposition, we will describe homomorphisms $\varphi:\widetilde{W}_n\to W_n$ explicitly. More specifically, we will decompose $\varphi$ into products of homomorphisms defined on $Q(R^{\vee})$, $N$, and $S_+$.  This will allow us to better understand these homomorphisms as we will describe how they impact the components of $\widetilde{C}_n$. 


    Following the convention of Bourbaki in {\rm\cite{bourbaki_groupes_1981}}, the coroot lattice of $W_n$ is
    $$ Q(R^{\vee})=\begin{cases} 
    \ZZ^n & W_n=C_n,\\
    \{(x_1,\dots,x_n)\in \ZZ^n:\sum x_i\equiv 0 \mod 2 \} & W_n=B_n,D_n.
    \end{cases}$$
    Define $\tau:\ZZ^n\to N$ by 
    $$ \tau(\lambda_1 e_1+\dots +\lambda_ne_n)=r_{1}^{\lambda_1}\dots r_n^{\lambda_n} $$
    this is an epimorphism on the coroot lattices of $B_n$, $C_n$, and $D_n$, the kernel of which is $(2\ZZ)^n$. Recall from Equation \eqref{equ: sgn:R-><g>} that, since $N$ is a subgroup generated by reflections, for any $w\in W_n$ with $w^2=1$, one can define $\mathrm{sgn}_N^{w}:\{\pm1\}\to \lb w\rb\leqslant  W_n$ to be the unique surjective homomorphism with kernel $N^+$. This map is trivial if and only if $w=1$. 
    In particular, because $z_0=r_1\dots r_n=-1$ is an involution, $\mathrm{sgn}_N^{z_0}:N\to N$ is the homomorphism  
    $$ \mathrm{sgn}_N^{z_0}(r_1^{\mu_1}\dots r_n^{\mu_n})=(r_1\dots r_n)^{\mu_1+\dots+\mu_n} .$$
    When $W_n=D_n$, $\mathrm{sgn}_N^{z_0}$ is the trivial map; but when $W_n=B_n$ or $C_n$, $\im(\mathrm{sgn}_N^{z_0})=\{\pm1\}=Z(W_n)\unlhd W_n$ and $\ker(\mathrm{sgn}_N^{z_0})=(W_n)_2^+$. Note that $\mathrm{sgn}_N^{z_0}\circ \tau$ is trivial on $B_n$, but nontrivial on $C_n$.
    We will combine $\tau$ and $\mathrm{sgn}_N^{z_0}$, and $\mathrm{sgn}_N^{w}$ to specify maps $\widetilde{W}_n\to W_n$, and Lemma \ref{lem: homs for nested semidirects} will ensure that these combinations are homomorphisms. 

    \begin{lemma}
        \label{lem: homs for nested semidirects}
        Let $A$, $B$ and $H$ be groups, with $A$ and $B$ abelian, and suppose $\alpha_1,\alpha_2:A\to B$ and $\beta_1,\beta_2:B\to B$ are homomorphisms. The map $f:A\rtimes (B\rtimes H)\to B\rtimes H$ defined by
        \begin{align*}
            (a,(b,h))\mapsto (\alpha_1(a)+\alpha_2(a)+\beta_1(b)+\beta_2(b),h)
        \end{align*}
        is a homomorphism if and only if the following conditions are satisfied:
        \begin{enumerate}[label=(\roman*)]
            \item $\al_i((b,h)\cdot a)=h\cdot \al_i(a)$ for all $a\in A$, $b\in B$, $h\in H$, and $i\in\{1,2\}$.
            \item $\be_i(h\cdot b)=h\cdot \beta_i(b)$ for all $b\in B$, $h\in H$, and $i\in\{1,2\}$.
        \end{enumerate}
    \end{lemma}

    For an arbitrary $\varphi \in \Hom(\widetilde{W}_n,W_n)$ and $(v,x,\sigma)\in \widetilde{W}_n$, Theorem \ref{thm: explicity homs types BCD} explicitly states  $\varphi(v,x,\sigma)$. The proof of this theorem involves straightforward checks that the stated maps are group homomorphisms using Lemma \ref{lem: homs for nested semidirects}, and that the kernels and images of the maps cover all possibilities listed in Table \ref{tab: counting Hom(W~,W) for A,B,C,D}. However, constructing these maps and developing a unified and simple way to describe them was complex.   

    \begin{theorem}
    \label{thm: explicity homs types BCD}
        Let $n\geq3$ and $n\neq 4$ and suppose $W_n=B_n, C_n$, or $D_n$ and $\varphi:\widetilde{W}_n\to W_n$ is a homomorphism. Then one of the following is true:
        \begin{enumerate}[label=(\alph*)]
            \item  The image-order of $\varphi$ satisfies the bound $|\im(\varphi)|\leq 8$ and there exist $w_1,w_2,w_3\in W_n$ satisfying $w_i^2=(w_iw_j)^2=1$ for all $1\leq i\leq j\leq3$ so that $\varphi$ has the form
            $$ \varphi(v,x,\sigma)=
            \mathrm{sgn}_{N}^{w_1}(\tau(v))\
            \mathrm{sgn}_{N}^{w_2}(x)\
            \mathrm{sgn}_{S_+}^{w_3}(\sigma). $$
            Moreover, the image and kernel of $\varphi$ depend on $(w_1,w_2,w_3)$ in the way described in Table \ref{tab: Image and kernel of small Hom(W~,W)}.
            
            \item The image-order of $\varphi$ satisfies the bound $n!\leq |\im(\varphi)|\leq 2^nn!$ and there exists $\epsilon_1,\epsilon_2,\epsilon_3,\epsilon_4\in \{0,1\}$ so that, up to postcomposition with an automorphism of $W_n$, $\varphi$ has the form
            $$ \varphi(v,x,\sigma)=
            \mathrm{sgn}_{N}^{z_0}(
            \tau(v)^{\epsilon_1}x^{\epsilon_2})\ 
            \tau(v)^{\epsilon_3}\ x^{\epsilon_4}\ \sigma .$$
            Moreover, the image and kernel of $\varphi$ depend on the tuple $(\epsilon_1,\epsilon_2,\epsilon_3,\epsilon_4)$ in the way described in Table \ref{tab: Image and kernel of large hom(W~,W)}. 
        \end{enumerate}
    \end{theorem} 
            
        \begin{table}[!ht]
            \centering
            \begin{tabular}{c  c  c}
             $(w_1,w_2,w_3)\in (C_n)^3$ & $\im(\varphi)$ & $\ker(\varphi)$ \\
        \noalign{\vskip 2pt}\hline\hline\noalign{\vskip 2pt}
            \multicolumn{3}{l}{$W_n=B_n$, $C_n$, or $D_n$:}\\ 
                $(1,1,1)$ & $\{1\}$ & $(W_n,Q(R^{\vee}),0)$  \\
                $(1,1,w_3)$ & $\FC_2$ & $\begin{cases}
                (\prescript{}{+}{W_n},Q(R^{\vee}),0) & W_n=B_n,C_n,\\
                (D_n^+,Q(R),0) & W_n=D_n.
                \end{cases}$ \\ 
        \noalign{\vskip 2pt}\hline\noalign{\vskip 2pt}
            \multicolumn{3}{l}{$W_n=B_n$ or $C_n$:}\\ 
                $(1,w_2,1)$ & $\FC_2$ & $((W_n)_1,Q(R^{\vee}),0)$ \\ 
                $(1,w_2,w_2)$ & $\FC_2$ & $(W_n^+,Q(R^{\vee}),0)$ \\ 
                $(1,w_2,w_3)$ & $\FC_2^2$ & $((W_n)_1^+,Q(R^{\vee}),0)$ \\ 
        \noalign{\vskip 2pt}\hline\noalign{\vskip 2pt}
            \multicolumn{3}{l}{$W_n=C_n$:}\\ 
                $(w_1,1,1)$ & $\FC_2$ & $(C_n,Q(R),(0,0))$ \\ 
                $(w_1,1,w_1)$ & $\FC_2$ & $(C_n,Q(R),(0,\al_n^{\vee}))$ \\ 
                $(w_1,w_1,1)$ & $\FC_2$ & $(C_n,Q(R),(\al_n^{\vee},0))$  \\ 
                $(w_1,w_1,w_1)$ & $\FC_2$ & $(C_n,Q(R),(\al_n^{\vee},\al_n^{\vee}))$  \\ 
                $(w_1,1,w_3)$ & $\FC_2^2$ & $(\prescript{}{+}{C_n},Q(R),0)$ \\ 
                $(w_1,w_2,1)$ & $\FC_2^2$ & $((C_n)_1,Q(R),0)$ \\
                $(w_1,w_1,w_3)$ & $\FC_2^2$ & $(\prescript{}{+}{C_n},Q(R),\al_n^{\vee})$ \\ 
                $(w_1,w_2,w_1)$ & $\FC_2^2$ & $((C_n)_1,Q(R),\al_n^{\vee})$\\
                $(w_1,w_2,w_2)$ & $\FC_2^2$ & $(C_n^+,Q(R),0)$ \\
                $(w_1,w_2,w_1w_2)$ & $\FC_2^2$ & $(C_n^+,Q(R),\al_n^{\vee})$ \\
                $(w_1,w_2,w_3)$ & $\FC_2^3$ & $((C_n)_1^+,Q(R),0)$ \\
            \end{tabular}
            \caption{The images and kernels of $\varphi:\widetilde{W}_n\to W_n$ for $W_n=B_n$, $C_n$, or $D_n$ when $|\im(\varphi)|\leq8.$ Note that $w_i\neq1$ and $w_i\neq w_j$ when $i\neq j$. If $W_n=D_n$, we can assume that $w_1=w_2=1$, and if $W_n=B_n$ we can assume that $w_1=1$, so for $D_n$ it suffices to consider only the first two rows, and for $B_n$ it suffices to only consider the first five rows.}
            \label{tab: Image and kernel of small Hom(W~,W)}
        \end{table}
        
        \begin{longtable}[h]{ c  l  l}
            $(\eps_1,\eps_2,\eps_3,\eps_4)$ & $\im(\varphi)$ & $\ker(\varphi)$ \\
        \noalign{\vskip 2pt}\hline\hline\noalign{\vskip 2pt}
            \multicolumn{3}{l}{$W_n=B_n$, $C_n$, or $D_n$:}\\ 

                $(0,0,0,0)$ & $\FS_n$ & $(N,Q(R^{\vee}),0)$ \\

                $(0,0,0,1) $ & $W_n$ & 
                $(\{1\},Q(R^{\vee}),0)$ \\ 
                
                $(0,0,1,0) $ & $W_n$ & $\begin{cases}((W_n)_2,2Q(R^{\vee}),0) & W_n=B_n,C_n,\\ 
                ((D_n)_{n-1,n},2Q(R)+2\mathbb{Z}\omega_1,0) & W_n=D_n.
                \end{cases}$ \\ 

                $(0,0,1,1) $& $W_n$ & $\begin{cases}
                ((B_n)_2^+,2Q(R),\alpha_{n-1}) & W_n=B_n, \\ 
                ((C_n)_2,2Q(R^{\vee}),\al_n^{\vee}) & W_n=C_n,\\[5pt]
                ((D_n)_{n-1,n},2Q(R)+2\mathbb{Z}\omega_1,\alpha_{n-1})& W_n=D_n.
                \end{cases}$  \\     
        \noalign{\vskip 2pt}\hline\noalign{\vskip 2pt}
            \multicolumn{3}{l}{$W_n=B_n$ or $C_n$:}\\ 
                
                $(0,1,0,0) $ & $\FC_2\times \FS_n$ &$((W_n)_2^+,Q(R^{\vee}),0)$ \\ 

                $(0,1,0,1) $  & 
                $\begin{cases}
                W_n/\{\pm1\}& \text{$n$ odd,}\\
                W_n& \text{$n$ even.}
                \end{cases}$ & 
                $\begin{cases}
                (\{\pm1\},Q(R^{\vee}),0)& \text{$n$ odd},\\
                (\{1\},Q(R^{\vee}),0)& \text{$n$ even.}
                \end{cases}$ \\[5pt]

                 $(0,1,1,0) $ & $
                 \begin{cases}
                    D_n & W_n=B_{n},\; n\text{ even},\\ 
                    \FC_2\times D_n & W_n=B_{n},\; n\text{ odd},\\
                     C_n & C_n.                         
                 \end{cases}$ & $\begin{cases}
                 ((B_n)_2^+,2Q(R),0) & W_n=B_{n},\; n\text{ odd},\\
                 ((B_n)_2,2Q(R),2\omega_n) & W_n=B_{n},\; n\text{ even},\\[5pt]
                 ((C_n)_2,2Q(R^{\vee}),2\omega_n^{\vee})& W_n=C_n.
                 \end{cases}$  \\ 

                $(0,1,1,1) $ & $
                 \begin{cases}
                    B_n/\{\pm1\} & W_n=B_{n},\; n\text{ odd},\\
                    B_n & W_n=B_{n},\; n\text{ even},\\
                    C_n & C_n.                         
                 \end{cases}$ & $\begin{cases} 
                 ((B_n)_2,2Q(R),2\omega_n+\alpha_n)& W_n=B_{n},\; n\text{ odd},\\
                 ((B_n)_2^+,2Q(R),\alpha_{n-1})& W_n=B_{n},\; n\text{ even},\\
                 ((C_n)_2,2Q(R^{\vee}),2\omega_n^{\vee}+\al_n^{\vee}) & W_n=C_n.
                 \end{cases}$  \\[5pt] 

        \noalign{\vskip 2pt}\hline\noalign{\vskip 2pt}
                \multicolumn{3}{l}{$W_n=C_n$:}\\ 
                
                $(1,0,0,0) $ & $\FC_2\times \FS_n$ & $((C_n)_2,Q(R),0)$  \\
                
                $(1,0,0,1) $ & $C_n$ & $(\{\pm1\},Q(R),\al_n^{\vee})$ \\ 
                                    
                $(1,0,1,0) $ & 
                $\begin{cases}
                C_n/\{\pm1\}& \text{$n$ odd,}\\
                C_n& \text{$n$ even.}
                \end{cases}$ & 
                $\begin{cases}
                ((C_n)_2,2P(R^{\vee}),0)& \text{$n$ odd,}\\
                ((C_n)_2,2Q(R^{\vee}),0)& \text{$n$ even.}
                \end{cases}$ \\[15pt]
                
                $(1,0,1,1) $ & $C_n$ & 
                $\begin{cases}
                ((C_n)_2^+,2P(R^{\vee}),2\alpha_{n-1})& \text{$n$ odd,}\\
                ((C_n)_2,2Q(R^{\vee}),2\omega_n^{\vee}+\al_n^{\vee})& \text{$n$ even.}
                \end{cases}$ \\[15pt]
                
                $(1,1,0,0) $ & $\FC_2\times \FS_n$  & $((C_n)_2,Q(R),2\al_{n})$  \\                     
                $(1,1,0,1) $ & $C_n$ & 
                $\begin{cases}
                (\{\pm 1\},Q(R),0) & \text{$n$ odd,}\\
                (\{\pm1\},Q(R),\al_n^{\vee}) & \text{$n$ even.} 
                \end{cases}$  \\[15pt]                                                           
                $(1,1,1,0) $ & $C_n$ & 
                $\begin{cases}
                ((C_n)_2^+,2P(R^{\vee}),0)& \text{$n$ odd},\\
                ((C_n)_2,2Q(R^{\vee}),2\omega_n^{\vee})& \text{$n$ even.}                    
                \end{cases}$  \\[15pt]

                $(1,1,1,1) $ & 
                $\begin{cases}
                C_n/\{\pm1\}& \text{$n$ odd,}\\
                C_n& \text{$n$ even.}
                \end{cases}$ & 
                $\begin{cases}
                ((C_n)_2,2P(R^{\vee}),2\al_{n-1})& \text{$n$ odd,}\\
                ((C_n)_2,2Q(R^{\vee}),\al_n^{\vee})& \text{$n$ even.}
                \end{cases}$
                \\ 
            \caption{The images and kernels of $\varphi:\widetilde{W}_n\to W_n$ when $n!\leq|\im(\varphi)|\leq 2^nn!$. Every $\varphi:\widetilde{W}_n\to W_n$ is the composition of a homomorphism $\rho_{(\eps_1,\dots,\eps_4)}:\widetilde{W}_n\to W_n$ corresponding to a row in this table with an automorphism $f\in \Aut(W_n)$. When $W_n=C_n$ and $n$ is odd, $\rho_{(\eps_1,\dots,\eps_4)}$ is unique and $f$ is unique except in the case where the image is isomorphic to $\FC_2\times\FS_n$. When $W_n=C_n$ and $n$ is even, there are six kernels that appears twice in the rows. Homomorphisms with these kernels are related by postcomposition by an outer automorphism of $C_n$. If $W_n=D_n$ then since $\mathrm{sgn}_N^{z_0}$ is trivial, $\eps_1$ and $\eps_2$ do nothing so it suffices to consider only the first four rows of the table. Similarly for $B_n$ we need only the first eight rows because $\mathrm{sgn}_N^{z_0}\circ\tau$ is trivial.}
            \label{tab: Image and kernel of large hom(W~,W)}
        \end{longtable}   

         \begin{proof}
        Suppose $W_n=C_n$. The abelianisation of $\widetilde{C}_n$ is isomorphic to $\FC_2^3$, where the three factors record the `sign' of elements in $Q(R^{\vee})$, $(C_n)_2$, and $S_+$. A small-image homomorphism is a composition of the abelianisation map with some homomorphism $\FC_2^3\to C_n$, and hence is a homomorphism. 
        
        We now focus on large-image homomorphisms. To show these maps are homomorphisms, it suffices to check that $\tau$ satisfies property $(i)$ and $\mathrm{sgn}_{N}^{z_0}$ satisfies property $(ii)$ of Lemma \ref{lem: homs for nested semidirects}. Let $v=\sum_{i=1}^n \lambda_ie_i \in Q(R^{\vee})$, $x=\prod_{i=1}^n r_i^{\mu_i}\in (C_n)_2$, and $\sigma\in S_+\cong \FS_n$, then 
        $$ 
        \tau((x,\sigma)\cdot v)=\tau\big((-1)^{\mu_1}\lambda_1e_{\sigma(1)}+\dots+(-1)^{\mu_n}\lambda_ne_{\sigma(n)}\big)=
        r_{\sigma(1)}^{\lambda_1}\dots r_{\sigma(n)}^{\lambda_n}=\sigma \cdot (r_1^{\lambda_1}\dots r_n^{\lambda_n})=\sigma\cdot \tau(v).
        $$
        Similarly, let $x=\prod_{i=1}^n r_i^{\mu_i}\in (C_n)_2$ and $\sigma\in \FS_n$ and observe that since $z_0=r_1\dots r_n$, the permutation $\sigma$ acts trivially on $\lb z_0\rb$, so we have
        $$ \mathrm{sgn}_{N}^{z_0}(\sigma\cdot (r_1^{\mu_1}\dots r_n^{\mu_n}))=\mathrm{sgn}_{N}^{z_0}\big(r_{\sigma(1)}^{\mu_1}\dots r_{\sigma(n)}^{\mu_n}\big)=z_0^{\mu_1+\dots+\mu_n}=\sigma\cdot z_0^{\mu_1+\dots+\mu_n}=\sigma\cdot \mathrm{sgn}_{N}^{z_0}(x) .$$
        
        It remains to compute the kernels of the homomorphisms corresponding each tuple $(\eps_1,\dots,\eps_4)$ and verify that these cover each of the possible kernels of homomorphisms $\widetilde{C}_n\to C_n$ listed in Table \ref{tab: counting Hom(W~,W) for A,B,C,D}. We will compute the kernel of the homomorphism $\varphi $ corresponding to $(0,0,1,1)$, the other cases are similar. 

        Suppose $(v,x,\sigma) \in \ker(\varphi)$, then 
        $$\varphi(v,x,\sigma)=(\tau(v)x,\sigma) =(1,1)\in \FC_2^2\rtimes \FS_n\cong C_n$$ 
        so $\sigma$ is trivial and $\tau(v)=x$. There exist $\mu_1,\dots,\mu_n\in \{0,1\}$ so that $x=r_1^{\mu_1}\dots r_n^{\mu_n}$,  and $\tau(v)=x$ if and only if 
        $ v=\sum_{i=1}^{n} (\lambda_i+\mu_i)e_i$ for some $\lambda_i\in 2\mathbb{Z}.$ 
        
        On the other hand, for $((C_n)_2,2Q(R^{\vee}),\al_n^{\vee})$ and since $c(r_n)=\al_n^{\vee}$, by Remark \ref{rmk: correspondence with triples},
        $$c(r_1^{\mu_1}\dots r_n^{\mu_n})=\mu_1e_1+\dots+\mu_n e_n.$$
        Hence, using Equation \eqref{equ: normal subgroup N~(K,L,c)}, we can write $((C_n)_2,2Q(R^{\vee}),\al_n^{\vee})$ as 
        \begin{align*} 
        ((C_n)_2,2Q(R^{\vee}),\al_n^{\vee})=&
        \big\{(v+ \mu_1e_1+\dots +\mu_ne_n,\: r_1^{\mu_n}\dots r_n^{\mu_n},\: 1): v\in 2Q(R^{\vee}), \mu_i\in\{0,1\} \big\}\\ 
        =&\big\{((\la_1+\mu_1)e_1+\dots +(\la_n+\mu_n)e_n,\: r_1^{\mu_n}\dots r_n^{\mu_n},\: 1): \la_i\in 2\ZZ, \mu_i\in\{0,1\} \big\}\\ 
        =& \ker(\varphi).
        \end{align*}

        The proofs for $B_n$ and $D_n$ are similar. However, recall that for $B_n$, $\mathrm{sgn}_{N}^{z_0}\circ \tau$ is trivial and for $D_n$, $\mathrm{sgn}_{N}^{z_0}$ is trivial. Hence for these groups we can assume $\eps_1=0$, and $\eps_1=\eps_2=0$, respectively.
    \end{proof}   

    We can use Theorem \ref{thm: explicity homs types BCD} to describe homomorphisms $\widetilde{W}_n\to W_n$ in terms of the standard generators $s_0,s_1,\dots,s_n$ of $\widetilde{W}_n$. Written as tuples $(v,x,\sigma)\in Q(R^{\vee})\rtimes (N\rtimes W_n)$, the standard generators are
    \begin{align*}
        s_0&= \begin{cases}
            (e_1+e_2,r_1r_2,(1,2)) & W_n=B_n,D_n,\\
            (e_1,r_1,1) & W_n=C_n.\\
        \end{cases}\\
        s_i&=(0,1,(i,i+1)) \qquad \text{for }i\in \{1,\dots,n-1\}.\\ 
        s_n&= \begin{cases}
            (0,r_n,1) & W_n=B_n,C_n,\\
            (0,r_{n-1}r_n,(n-1,n)) & W_n=D_n.\\
        \end{cases}\\
    \end{align*}
    The proof of Corollary \ref{cor: Hom(Wn~,Wn) standard gens} amounts to evaluating the homomorphisms defined in Theorem \ref{thm: explicity homs types BCD} on $s_0,\dots, s_n$. For example, when $W_n=C_n$ we see that
    $$ \varphi_{(0,1,1,0)}(s_0)=z_0r_{\widetilde{\alpha}} \qquad 
        \varphi_{(0,1,1,0)}(s_i)=s_i \quad i\in \{1,\dots,n-1\} \qquad 
        \varphi_{(0,1,1,0)}(s_n)=z_0. $$
    
    \begin{corollary}[Homomorphisms in standard generators]
        \label{cor: Hom(Wn~,Wn) standard gens}
        Suppose $W_n=B_n$, $C_n$, or $D_n$. Let $n\geq3$ and $n\neq 4$, and $\varphi\in \Hom(\widetilde{W}_n,W_n)$.
        \begin{enumerate}[label=(\alph*)]
        \item If $W_n=B_n$, then one of the following is true.
            \begin{enumerate}[label=(\roman*)]
                \item There exist commuting $w_1,w_2\in B_n$ satisfying $w_i^2=1$ so that $\varphi$ is:
                $$  \varphi(s_i)=w_1 \quad \text{for }i=0,\dots,n-1; \quad \varphi(s_n)=w_2.$$
                \item There exist $\eps_1,\eps_2,\eps_3\in\{0,1\}$ so that, up to postcomposition with an automorphism of $B_n$, $\varphi$ is:
                $$
                \varphi(s_0)=r_{\widetilde{\alpha}}^{\eps_1}s_1^{\eps_1+1} \quad \varphi(s_i)=s_i \quad \text{for }i=1,\dots,n-1; \quad \varphi(s_n)=z_0^{\eps_2}s_n^{\eps_3}.
                $$
            \end{enumerate}
        \item If $W_n=C_n$, then one of the following is true.
            \begin{enumerate}[label=(\roman*)]
                \item There exist pairwise commuting $w_1,w_2,w_3\in C_n$ satisfying $w_i^2=1$ so that $\varphi$ is:
                $$ \varphi(s_0)=w_1\quad \varphi(s_i)=w_2 \quad \text{for }i=1,\dots,n-1; \quad \varphi(s_n)=w_3.$$
                \item There exist $\eps_1,\eps_2,\eps_3,\eps_4\in\{0,1\}$ so that, up to postcomposition with an automorphism of $C_n$, $\varphi$ is:
                $$
                \varphi(s_0)=z_0^{\eps_1}r_{\widetilde{\alpha}}^{\eps_2} \quad \varphi(s_i)=s_i \quad \text{for }i=1,\dots,n-1; \quad \varphi(s_n)=z_0^{\eps_3}s_n^{\eps_4}.
                $$
            \end{enumerate}
        \item If $W_n=D_n$, then one of the following is true.
            \begin{enumerate}[label=(\roman*)]
                \item There exists $w\in D_n$ with $w^2=1$ so that $\varphi$ is:
                $$  \varphi(s_i)=w \quad \text{for }i=0,\dots,n.$$
                \item There exist $\eps_1,\eps_2\in\{0,1\}$ so that, up to postcomposition with an automorphism of $D_n$, $\varphi$ is:
                $$
                \varphi(s_0)=r_{\widetilde{\alpha}}^{\eps_1}s_1^{\eps_1+1} \quad \varphi(s_i)=s_i \quad \text{for }i=1,\dots,n-1; \quad \varphi(s_n)=s_{n-1}^{\eps_2}s_n^{\eps_2+1}.
                $$
            \end{enumerate}
        \end{enumerate}
    \end{corollary}

\printbibliography[title={References}]

\end{document}